\documentclass[a4paper,10pt]{article}
\usepackage{amsfonts}
\usepackage{graphicx,amssymb,amsmath}
\usepackage{amscd}
\usepackage[english]{babel}
\usepackage[ansinew]{inputenc}
\usepackage{slashbox}
\usepackage[all]{xy}
\usepackage{latexsym}
\usepackage{amsmath}
\usepackage{stmaryrd}

 \hsize \textwidth
 \advance \hsize by -\marginparwidth
\newcommand{\field}[1]{\mathbb{#1}}
\newcommand{\R}{\field{R}}
\newcommand{\Z}{\field{Z}}

\newcommand{\N}{\field{N}}
\newcommand{\M}{\field{M}}

\newcommand{\E}{\field{E}}

\newtheorem{theo}{Theorem}[section]

\newtheorem{defi}[theo]{Definition}
\newtheorem{lem}[theo]{Lemma}

\newtheorem{prop}[theo]{Proposition}
\newtheorem{rem}[theo]{Remark}

\newtheorem{conv}[theo]{Convention}
\newtheorem{exe}[theo]{Example}
\newtheorem{exes}[theo]{Examples}

\def\dis{\ds}
\def\ap{\rightarrow}
\def\S{\Sigma}
\def\l{\lambda}
 \def\a{\alpha}
\def\b{\beta}
\def\g{\gamma}
\def\G{\Gamma}
\def\L{\Lambda}

\def\t{\tau}
\def\d{\delta}
\def\o{\omega}

\def\l{\lambda}
\def\n{\nu}
\def\f{\flat}

\def\r{\rho}

\def\m{\mu}
\def\n{\nu}
\def\s{\sigma}

\def\so{\underline}

\def\O{\Omega}

\def\e{\epsilon}

\def\~{\tilde}
\def\dis{\displaystyle}
\input epsf

\title{{Almost Lie  structures on an anchored Banach bundle}}
 \author{ P Cabau \& F. Pelletier}

\date{}

\begin{document}

\maketitle

\begin{abstract}
 Under appropriate assumptions, we generalize the concept of linear almost Poisson structures, almost Lie algebroids,  almost differentials in the framework of Banach anchored bundles and the relation between these objects. We then obtain an adapted formalism for mechanical systems which is illustrated by the evolutionary problem of the "Hilbert snake" as exposed in \cite{PeSa}.
  \end{abstract}
\section{Introduction}
Recent developments about geometric formalism on anchor bundles on a finite dimensional manifold have helped to build a general framework for studying mechanical systems. Essentially, these geometric structures concern,  linear almost Poisson structures, almost Lie algebroids and almost differentials (see for example \cite{GLMM}, \cite{GLMM}, \cite{LMM}, \cite{Marl}, \cite{Marr}, \cite{PoPo} and all references inside these papers).\\

The purpose of this paper is to give a generalization of these geometric structures in the context of Banach anchored bundles. Of course, this framework leads to a lot of obstructions. At first, any local section of a Banach bundle cannot be extended to a global section without some properties of regularity of the typical fiber. So, in the setting of "Banach algebroid", we must impose that the Lie Bracket of sections (defined globally) has a property of "localization" . Indeed, if the Banach manifold is regular, this property is always satisfied, as in finite dimension. However, in the general case, we must impose such a condition (see section \ref{loc}). On the other hand hand, for a Poisson structure on a Banach manifold, we meet the same type of problem but also, if the typical Banach model is not reflexive, we must impose some other conditions (see section \ref{APmor}). Finally, the most important obstruction appears in the context of "Lie differential":  not only we still meet the problem of localization of sections but, on the opposite of finite dimension, the graded algebra of forms on a Banach space is not generated by elements of degree zero and degree one, so, in general, such a differential is not characterized by its values on elements of these types.\\

However, by imposing appropriate assumptions, we defined the concept of almost Lie bracket, almost Lie algebroid which is a generalization of Lie algebroid on Banach manifold introduced in \cite{Ana} and \cite{Pel}. Now, recall that in finite dimension, on one hand  there exists a bijection between Lie algebroid structures on an anchored bundle and  Poisson structures on its dual, and a bijection between Lie algebroid structures  and Lie differentials (see for instance \cite{Marl} or \cite{GLMM} among many references). In fact, in our context, if the typical fiber of the anchored bundle is not reflexive, we   have a bijection with almost Lie algebroid
on an anchored bundle  but here with "sub almost Poisson structure" on its dual, that is, such a structure is only defined on the set of sections of a "canonical subbundle" of the dual bundle (see subsection \ref{alg-Poisson}). Of course, in the framework of paracompact Hilbert manifolds, all these obstructions do not exist and we recover the general setting of the finite dimensional context. On the opposite, in the infinite dimensional context, we do not have  any bijection between almost Lie algebroid structures and almost Lie differentials even under appropriate assumptions (see subsection \ref{A-deriv}).\\

In the  following section, we recall the concept of graded exterior algebra, some classical properties of Banach manifolds under which the mentioned previous assumptions  can be avoid. We also precise our notations in local context. The notion of almost Lie algebroid  (AL algebroid in short) on a Banach anchored bundle is developed in section \ref{ALalg}. The section \ref{APmorp} is devoted to the relation between "sub almost Poisson structure" and AL algebroid in one hand and almost Lie differential and AL algebroid on the other hand. At first we expose the context of sub  almost Poisson morphism (sub AP morphism in short) (subsection \ref{APmor}). Then we look for the equivalence of AP morphism and AL algebroid structure (subsection \ref{alg-Poisson}). Finally, under strong appropriate assumption on an almost differential, we can associate an AL algebroid structure (subsection \ref{A-deriv}).
As applications of our results,  in section \ref{mecha}, we look for an adaptation of  classical formalism for mechanical systems:  Hamiltonian system, Hamilton-Jacobi equation, Lagrangian system and Euler-Lagrange equation. Finally, in section \ref{exemple}, we illustrate this formalism in the context of the evolution of the "head" of  a "Hilbert snake"  form the results of \cite{PeSa}.
\section{Preliminaries and notations}\label{prelnot}

\subsection{Graded exterior algebra on a Banach space}

Given a  Banach space $\E$, we denote by  $\L^k\E^*$ the Banach  space of exterior forms of order $k$ of $\E$.  More precisely,  the set $\L^k\E^*$  can   be identified with the closed subspace  of $k$-multilinear skew-symmetric maps in the Banach space ${\cal L}^k(E)$ of $k$-multilinear forms on $\E$. This space is the closure of the vector space generated by all exterior products of $1$ forms  $\{\xi_{i_1}\wedge\cdots\wedge\xi_{i_k},\; i_1<\cdots<i_k\}$.  (for a complete description, see \cite{Ram}).  Set
$$\L\E^*=l^{\infty}(\dis\bigoplus_{k=0}^\infty \L^k\E^*)=\{\o=\dis\sum_{k=0}^\infty \o_k \textrm{ with } \o_k\in \L^k\E^*,\; \dis\sup_{k}||\o_k||<\infty\}.$$
Then, $ \L\E^* $ is a Banach space which, provided with the exterior product,  is an  algebra .\\

If we consider the Banach space $\E$, isometrically embedded in $\E^{**}$, we can define, in the same way, the vector space  $\L^k\E$  spanned by all exterior products  $\{u_{i_1}\wedge\cdots\wedge u_{i_k},\; i_1<\cdots<i_k\}$.  So, if we set:
$$\L\E=l^1(\dis\bigoplus_{k=0}^\infty \L^k\E)=\{u=\dis\sum_{k=0}^\infty u_k  \textrm{ with } u_k\in \L^k\E,\;\dis\sum_{k=0}^\infty||u_k||<\infty\}$$
then,  $\L\E$ is a Banach space, which, provided with the exterior product is a Banach algebra and also a graded algebra. Moreover $\L^k\E^*$ is isomorphic to $(\L^k\E)^*$ for any $k\geq 0$ (see \cite{Ram}).\\
Notice that  we have $\Lambda^{0}\E^*=\L^0\E=\mathbb{R}$,  $\Lambda^{1}\mathbb{E}=\mathbb{E}$, and $\L^1\E^*=\E^*$.

\bigskip
\noindent The {\it  interior product}  of $\o\in
\Lambda^{k}\mathbb{E}^{\ast}$ by a vector  $v\in \E,$ denoted by $i_{v},$ is characterized, as usual in the following way:

\begin{description}
\item -- if $k\leq0$, $i_{v}=0$ on $\Lambda^{k}\mathbb{E}^{\ast}$

\item -- if $k=1$,
 then $i_{v}\o=\left\langle \o,v\right\rangle \in\mathbb{R}$

\item -- if $k>1$, for $v_{1}\dots v_{k-1}\in\mathbb{E}$, then
$i_{v}\o(  v_{1},\dots,v_{k-1})  =\o\left(  v,v_{1}
,\dots,v_{k-1}\right)$

and then $i_{v}\o\in\Lambda^{k-1}\mathbb{E}^{\ast}$.
\end{description}

 In fact,  given any fixed $v\in \mathbb{E}$,  the interior product  $i_v$ can be clearly  extended to a continuous endomorphism of  $\L\mathbb{E}^{\ast}$ which is a derivation of degree  $-1$  of $\L\mathbb{E}^{\ast}$  that is $ i_v$ is  send  each factor   $\Lambda^{k}\mathbb{E}^{\ast}$ of the graduation into  the factor $\Lambda^{k-1}\mathbb{E}^{\ast}$
 $\L\mathbb{E}^{\ast}$.

\noindent The interior product  $i_{P}$ by a multivector $P\in\Lambda
^{p}\mathbb{E}$ is defined in the following way:
\begin{description}
\item -- if  $k<0$, $i_{P}=0$

\item -- if $k=0$, $P\in \R$ and for any form $\alpha\in
\Lambda\mathbb{E}^{\ast}$, then
$i_{P}\o=P\o$

\item -- if  $k\geq1$ and if $P$ is decomposable, i.e.  $P=v_{1}\wedge\cdots\wedge v_{k}$, we set
$i_{v_{1}\wedge\cdots\wedge v_{k}}=i_{v_{1}}\circ\cdots\circ i_{v_{k}}$
\end{description}

\noindent  We can  extend, by linearity and continuity,  the definition of  $i_{P}$, for any
 $P\in\Lambda^{k}\mathbb{E}$,  to the graded algebra $\Lambda\mathbb{E}^{\ast}$. For any fixed $P\in\Lambda^{k}\mathbb{E}$ we then get 

 an endomorphism  $i_p$ of  degree $-p$ of the graded algebra
 $\Lambda\mathbb{E}^{\ast}$.\\

\noindent  Let  $\t:E\rightarrow M$ be a Banach bundle of typical fiber $\mathbb{E}$.
In this situation, we denote by:

-- ${\cal F}(E)$ the algebra of smooth functions on $E$;

-- $\G(\t)$ the $\mathcal{F} $--module of smooth sections $C^{\infty}$ of this bundle;

-- $\L^k\G^*(\t)$  the  $\mathcal{F}$--module of  sections of the Banach bundle $\L^kE^*$ of typical fiber $\L^k\E^*$;

-- $\L^k\G(\t)$ the  $\mathcal{F}$--module of sections of the Banach bundle $\L^kE$ of typical fiber $\L^k\E$;

-- $\L\G^*(\t)$  the  $\mathcal{F}$--module of  sections of the Banach bundle $\L E^*$ of typical fiber $\L \E^*$;

-- $\L\G(\t)$ the  $\mathcal{F}$--module of  sections of the Banach bundle $\L E$ of typical fiber $\L \E$.\\

\noindent For any $P$ in $\Lambda^{k}\G(\t)$, the integer  $k$ is called
the degree of  $P$ and we set $k=\deg P$.

\noindent Notice that  $\Lambda^{0}\G(\t)=\mathcal{F}$ and $\Lambda^{1}\G(\t)=\G(\t)$. For the exterior product of vectors (resp. of forms) we get  a structure of  graded exterior algebra on  $\L\G(\t)$ (resp.  $\L\G^*(\t)$)

\subsection{Some classical properties of   Banach manifolds}
First we recall some classical  properties of Banach manifolds. The reader can find complete references about all these properties in \cite{KrMi}. \\

 Let $M$ be a smooth Banach manifold  modeled on the Banach space $\M$. The manifold $M$ is {\bf paracompact} if and only if  the topology of $M$ is metrizable. In particular the Banach space $\M$ must be also paracompact. This is always true for any (eventually non separable) Hilbert space.
A Banach space which has  $C^k$ partitions of unity is called {\bf $\bf C^k$-paracompact}, for $k\in \N\cup{\infty}$. Any paracompact manifold modeled on a $C^k$-paracompact Banach space has also $C^k$-partitions of unity and so is $ C^k$-paracompact. \\

The Banach manifold $M$ is said {\bf  $C^k$-regular} (resp. {\bf smooth regular}) if for any $x\in M$, there exists an open neighborhood $U$ of $x$ and  a $C^k$ (resp. smooth) function $f:U\ap \R$  such that $f(x)=1$ and the closure of the set $\{z\;, f(z)\not=0\}$ is contained in $U$; such a function is called a {\bf bump function}.  Notice that $M$ is smooth-regular if and only if $\M$ is $C^k$-regular for any  $k\in \N\cup\infty$. \\
Of course not all Banach spaces are  $C^k$-regular for $k>0$: for instance, $l^1(\G)$, for any set $\G$, is not $C^1$ regular (see \cite {KrMi}).  But any $C^k$-regular Banach space is paracompact (for more details about regularity and paracompactness see also \cite{Arn}, \cite{Kun}, \cite{Llo},  \cite{Vand} and \cite{Ion}).\\

Notice that if $M$ is smooth  paracompact, then $M$ is smooth regular.

\subsection{Local coordinates in a Banach bundle}\label{loctriv}
Consider a Banach bundle $(E,\t,M)$ and the associated dual bundle $(E^*,\t_*,M)$. Fix $x\in M$  and consider any open neighborhood  $U$  of $x$  such that  $E_U$ is isomorphic to the trivial bundle $U\times\E$, which we always write $E_U\equiv U\times \E$; we then say that $E_U$ is trivialized. Then we also have

 $E_U^*\equiv U\times E^*$;

 $T^*E^*_{| E^*_U}\equiv U\times\E^{*}\times\M^*\times \E^{*}$;

 $TE^*_{| E^*_U}\equiv U\times\E^{*}\times\M\times \E^*$.  \\

Here we consider $\E$ as a Banach subspace of $\E^{**}$. Taking into account these equivalences, we get the following coordinates:

$s=(x,u)$ on $E_U\equiv U\times \E$

$\s=(x,\xi)$ on  $E^*_U\equiv U\times \E^*$

$(s,v)=(s,v_1,v_2)$ on  $TE_{| E_U}\equiv U\times\E\times\M\times \E$

$(\s,w)=(\s,w_1,w_2)$ on $TE^*_{| E^*_U}\equiv U\times\E^{*}\times\M\times \E^*$

$(\s,\eta)=(\s,\eta_1,\eta_2)$ on $T^*E^*_{| E^*_U}\equiv U\times\E^{*}\times\M^*\times \E^{**}$

\subsection{Local coordinates in basis}\label{locbas}

Suppose that the Banach space $\M$ has an unconditional, eventually uncountable, basis, $\{\m_i\}_{i\in I}$ and denote by $\{\m^*_i\}_{i\in I}$ the associated weak-$*$ basis of the dual $\M^*$ (see \cite{FiWo}). So, each $z\in \M$ can be written in a unique way:
$$x=\dis\sum_{i\in I}x^i \m_i.$$
Note that we have $x^i=\m^*_i(x)$.\\
On the other hand, each $\o\in \M^*$ can be "weak-$*$" written in a unique way:
$$\o\equiv\dis\sum_{i\in I}\o_i\m^*_i$$
which means that, for any $u\in \M$, we have:
$$<\dis\sum_{i\in I} \o_i \m^*_i,u>=<\o,u>$$
In fact we have $\o_i=<\o,\m_i>$.\\

Consider a chart $(U,\phi)$ on $M$. Via the diffeomorphism $\phi$, we can identify $U$ with an open set of $\M$ and  so any $x\in U$ can be written in a unique way
$x=\dis\sum_{i\in I} \m^*_i(x) \m_i$. We will say that the set of  maps $\{x^i:=\m^*_i:U\ap \R\}_{i\in I}$ is  the local system of coordinates on $U$.
As the tangent bundle $TM_{| U}$ is isomorphic to $U\times \M$, we denote by $\{\dis\frac{\partial}{\partial x^i}\}_{i\in I}$ the  basis of each fiber $T_zM$, for $z\in U$, canonically associated to $\{\m_i\}_{i\in I}$. So any vector field $X$ on $U$ can be written in a unique way as:
$$X=\dis\sum_{i\in I} X_i \dis\frac{\partial}{\partial x_i}$$
Moreover, as $X$ can be identified with a map from $U$ to $\M$ we have $X_i=\m^*_i\circ X$ and so each component $X_i$ is a smooth function.

In the same way, the cotangent bundle $T^*M_{| U}$ is isomorphic to $U\times \M^*$. We denote by $\{dx_i\}_{i\in I}$ the  weak-$*$ basis on each fiber $T_z^*M$, for $z\in U$, canonically associated to $\{\m^*_i\}_{i\in I}$. Again each $1$-form $\o$ on $U$ can be weak-$*$ written
$$\o\equiv\dis\sum_{i\in I}\o_idx_i$$
where of course we have
$$<\o,X>=<\dis\sum_{i\in I} \o_idx_i,X>.$$ Again, each component $\o_i$ is a smooth function.\\

On the other hand, consider a Banach bundle $\t:E\ap M$, and suppose that there exists an unconditional basis $\{e_\a\}_{\a\in A}$ for $\E$. Consider an open set $U\subset M$ which is a chart domain and such that $E_U\equiv U\times \E$.

With the previous properties, we denote again by $e_\a$ the constant section $x \mapsto e_\a$ in $E_U$. Each section  $s\in \G(\t_U)$ can be written as:
$$s=\dis\sum_{\a\in A}e^*_\a(s) e_\a$$
Again each "component" $u_\a=e^*_\a(s)$ is a smooth function on $U$.

So we have the following local coordinates on $E_U$:

$\bullet\;\;$ $(x,u)=(x^i,u^\a)$ on $E_U$

if $\M$ has also an unconditional basis, on $E_U$ the tangent space $T_sE_U$ is spanned by the basis $\{\dis\frac{\partial}{\partial x^i}\}_{i\in \N}$ and
$\{\dis\frac{\partial}{\partial u^\a}\}_{\a\in A}$  where $\{\dis\frac{\partial}{\partial u^\a}\}_{\a\in A}$ is identified with the basis $\{e_\a\}_{\a\in A}$ of $\{s\}\times \E$

So we have the following local coordinates on $TE_U$

$\bullet\;\;$ $(s,v_1,v_2)=(x^i,u^\a,X^i,U^\a)$

In the same way, for the dual bundle $\t_*:E^*\ap M$, on $E^*_U$ we have the constant sections $e^*_\a:x\mapsto e^*_\a$ and any  section $\s$  of $E^*_U$, we also can write, in a "weak-$*$" way,
$$\s=_{w}\dis\sum_{\a\in A}\xi_\a e^*_\a$$
and again each component  $\xi_\a$ is a smooth function on $U$.

So we have the following local coordinates on $E^*_U$

$\bullet\;\;$   $\s=(x,\xi)=(x^i,\xi_\a)$   ("weak-$*$ coordinates" for $\xi_\a$ )

if $\M$ also has an unconditional basis, on $E^*_U$,  the  tangent space  $T_\s E^*_U$ is spanned by the basis $\{\dis\frac{\partial}{\partial x^i} \}_{i\in I}$ and "weakly-$*$"

spanned by the basis  $\{ \dis \frac{\partial}{\partial \xi_\a}\}_{\a\in A}$ where again $\{ \dis \frac{\partial}{\partial \xi_\a}\}_{\a\in A}$ is identified with the weak-$*$ basis $\{e^*_\a\}_{\a\in A}$;

So we have the following local coordinates on $TE^*_U$

$\bullet\;\;$  $(\s,w_1,w_2)=(x^i,\xi_\a,X^i,\Xi_\a)$

\subsection{Derivations and vector fields}\label{deriv}
Let $M$ be a Banach manifold.
 Recall that a (global) {\bf derivation} of $\cal F$ is a $\R$-linear map $\partial :{\cal F}\ap {\cal F}$ such that:
$$\partial(fg)=f \partial(g)+ \partial(f) g$$
We denote by $\cal D$ the vector space of all derivations of $\cal F$.\\

An {\bf operational vector field} $\partial$ at $x\in M$  is a derivation of ${\cal F}(U)$ for some neighborhood $U$ of $x$  which is compatible with restriction to open $V\subset U$ i.e. $\partial$ induces a unique derivation $\partial_V$ of ${\cal F}(V)$ such that
$$\partial(f)_{| V}=\partial_V(f_{|V})$$

Let  $D_xM$  be the vector space of operational vector field at $x$ and $DM=\dis\cup_{x\in M} D_xM $ the set of operational vector fields.. In fact, the canonical projection $\hat{p}_M : DM \ap M$ gives rise to a structure of Banach bundle (see \cite{KrMi}). Unlike to the context of finite dimensional manifolds, if $M$ is not smooth regular, then the set of germs at $x\in M$ of elements of $\cal D$ can be smaller than ${\cal D}_x$. On the other hand, any local vector field on $M$ gives rise to an operational vector field but there exist elements of  ${\cal D}_x$ which do not induce local vector fields (for more details see \cite{KrMi}); in particular we have $TM \varsubsetneq DM$.

\section{Almost Banach Lie Algebroid }\label{ALalg}

\subsection{ Almost Lie bracket on an anchored Banach bundle }\label{ALbracket}

Let  $\t:E\rightarrow M$ be a Banach bundle of typical fiber $\mathbb{E}$. We will denote by $E_x=\t^{-1}(x)$ the fiber over $x\in M$.

A Banach morphism bundle $\rho:E\rightarrow TM$ is called an {\bf anchor}. This morphism induces a map, again denoted by $\r$
from $ \G(\t)$ to $\G(M) $ defined  for any $x\in M$ and any section $s$ of  $E$ by: $\r(s)(x)=\r\circ s(x)$. We say that $(E,\t,M,\rho)$ is an {\bf anchored Banach bundle}.\\

\bigskip
\bigskip
\bigskip

\noindent\so{\it  Local expressions} :${}$\\

In the context of local trivializations (see subsection \ref{loctriv}), we have:

\begin{eqnarray}\label{loctrivr}
\r(x,u)\equiv (x,u) \mapsto (x,R_x(u))
\end{eqnarray}
where $R: U\ap L(\E,\M)$.\\
  Suppose that the Banach spaces $\M$ and $\E$ have basis. According to subsection \ref{locbas},  on  any appropriate open $U$ in $M$,  any anchor $\r$ is locally characterized by a family $\{\r_\a^i\}_{i\in I,\a\in A}$ of smooth functions such that:
\begin{eqnarray}\label{locrho}
\r(e_\a)=\dis\sum_{i\in I}\r_\a^i\frac{\partial}{\partial x^i}
\end{eqnarray}

\bigskip


\begin{defi}${}$
\begin{enumerate}
\item  An almost Lie bracket (AL-bracket for short) on an anchored bundle  $\left(  E,\t,M,\rho\right)  $
is a  bracket $[  .,.]_\r  $ which satisfies the Leibniz property:
$$[ s_{1},fs_{2}]_\r  =f [  s_{1},s_{2}]_\r + \left(  {\rho}\left(  s_{1}\right)  \right)  \left(  f\right)\ s_{2}$$
for any $f\in\mathcal{F}$ and  $s_{1},s_{2}\in \G(\t)$.\\
In this situation, $(E,\t,M,\r,[.,.]_\r)$ is called an {\bf almost Lie Banach algebroid } (AL-algebroid for short).
\item  A Lie bracket (L-bracket for short) on an anchored bundle  $\left(  E,\t,M,\rho\right)  $
is an AL- bracket $[  .,.]_\r  $ which satisfies the Jacobi identity:
for all $s_1,s_2,s_3\in \G(\t),$
$$J(s_1,s_2,s_3)=[  s_1,[ s_2,s_3]  ]  +[ s_2,[ s_3,s_1]  ]  +[s_3,[ s_1,s_2] ]=0 $$
 In this case $\left(  E,\t,M,\rho,[  .,.]_\r  \right)  $ is called  a {\bf Lie Banach algebroid}  (L-algebroid  for short)
\end{enumerate}
\end{defi}

When $\left(  E,\t,M,\rho,[  .,.]_\r  \right)  $   is a L-algebroid, then the bracket  $[  .,.]_\r$ induces on  $\G(\t)$ a Lie algebra structure.
In this case
$\r:(\G(\t),[.,.]_\r) \ap (\G(M),[.,.])$
is a Lie algebra morphism.

\begin{exes}\label{ex1}${}$
\begin{enumerate}
\item Let $(E,\t,M)$ be a Banach subbundle of $(TM,p_{M},M)$ which is complemented, i.e. there exists a Banach subbundle $(F,p,M)$ of $(TM,p_{M},M)$ such that $T_xM=E_x\oplus F_x$. Let $\pi_1:TM\ap E$ be the Banach morphism associated to the projection of $T_xM$ onto $E_x$ whose kernel is $F_x$. We define
$$\left[  X,Y\right]  _{E}= {\pi_1}[X,Y]$$
where $[.,.]$ is the usual Lie bracket on vector fields.
 Then $\left(  E,\t,M,\rho,\left[  .,.\right]  _{E}\right)  $ is an AL-algebroid  where the anchor is the natural inclusion $\rho$ of $E$ in $TM$.
 Notice that $\left(  E,\t,M,\rho,\left[  .,.\right]  _{E}\right)$ is a L-algebroid if and only if  $(E,\t,M)$ is involutive. In particular $(TM,p_{M},M,Id,[.,.])$ is a L-algebroid.

 
 \item  Consider a smooth right action $\psi:M \times G\ap M$ of a connected  Banach Lie group $G$ over a Banach manifold $M$. Denote  by $\cal G$ the Lie algebra of $G$. We have a natural morphism $\rho$ from the trivial Banach bundle $M \times {\cal G}$ into $TM$ which is defined  by
$$\rho(x,X)=T_{(x,e)}\psi(0,X)$$
For any $X$ and $Y$ in $\cal G$, we have:
$$\rho(\{X,Y\})=[\rho(X),\rho(Y)]$$
  where $\{\;,\;\}$ denotes the Lie algebra bracket on $\cal G$ ( \cite{Bo}, \cite{KrMi}).
 On the trivial bundle $M\times {\cal G}$, each section can be identified with a map $\s: M\ap {\cal G}$  we define a Lie bracket on the set of sections by
 $$\{\{\s,\s'\}\}(x)=\{\s(x),\s'(x)\}+d\s(\xi_{\s'(x)})-d\s'(\xi_{\s(x)})$$
We get   an anchor $\Psi: M\times {\cal G}\ap TM$ by $\Psi(x,X)=\xi_X(x)$
  It follows that  $(M\times {\cal G},\Psi,M,\{\{\;,\;\}\})$ has  a  Banach  Lie algebroid structure  on $M$.
  
\item Let $\pi: N\ap M$ be a submersion between Banach manifolds. The subspaces $V_uN = T_u\pi^{-1}(x)\subset T_{x,u}N$,  denoted by $VN$ defined a Banach sub-bundle of $p_N:TN\ap N$   called the vertical subbundle. As the Lie bracket of two vertical vector fields is again a vertical vector field, we get a L-algebroid on $(VE, {\t_{E}}_{| VE},E)$.
    
\item Let $\theta$ be a $1$-form on a Banach manifold $M$ such that $d\theta$ is a weak symplectic form i.e.  the canonical map $\theta^\f: TM \ap T^*M$ defined by $\theta^\f(X)=i_X d\theta$ for any $X\in T_xM$ is injective (see \cite{OdRa2}). Assume that  $\theta^\f$ is closed  i.e.  $T^\f_x M=\theta^\f(T_xM)$ is closed in $T_x^*M$.
If $i^\f:T^\f M\ap T^*M$ is the natural inclusion, then we set $q^\f_M=q_M\circ i^\f:T^\f M\ap M$ the restriction of   $q_M:T^*M\ap M$. Then  $(T^\f M,q^\f_M,M)$ is a Banach subbundle of $(T^*M,q_M,M)$. Denote by $\Pi:T^\f M\ap TM$ the morphism $(\theta^\f)^{-1}$. We define  a structure of L-algebroid on the anchored bundle $(T^\f M,q^\f_M,M,\Pi)$ by setting
$$[\eta,\zeta]^\f=\theta^\f ([\Pi\eta,\Pi\zeta])$$
where $[.,.]$ is the usual Lie bracket of the  vector fields $\Pi\eta$ and $\Pi\zeta$. So $(T^*M,q_M,M,[.,.]^\f)$ is a L-algebroid.\\
This situation precisely occurs on the cotangent bundle $T^*M$ of any Banach manifold $M$ where $\theta$ is the Liouville $1$-form on $T^*M$ (see \cite{Lan})
\end{enumerate}
\end{exes}

\begin{defi}\label{morpalg}${}$\\
Let $(E_i,\t_i,M,\r_i,[.,.]_{\r_i})$, $\; i=1,2$, be two AL-algebroids (resp. L-algebroids). A morphism $\Psi$ from $(E_1,\t_1,M)$ to $(E_2,\t_2,M)$ (over $Id_M$) is called an {\bf AL-algebroid morphism} (resp. {\bf L-algebroid morphism}  if we have:
\begin{enumerate}
  \item $\Psi\circ \r_2=\r_1$
  \item $[\Psi(s_1),\Psi(s_2)]_{\r_2}=\Psi([s_1,s_2]_{\r_1})$ for any $s_1, s_2\in \G(\t_1)$
\end{enumerate}
\end{defi}

Notice that if $(E_i,\t_i,M,\r_i,[.,.]_{\r_i})$, $\; i=1,2$, are two L-algebroids, any AL-algebroid morphism $\Psi$ from $(E_1,\t_1,M,\r_1,[.,.]_{\r_1})$ to $(E_2,\t_2,M,\r_2,[.,.]_{\r_2})$  induces a Lie algebra morphism from $(\G(\t_1), [.,.]_{\r_1})$ to $(\G(\t_2), [.,.]_{\r_2})$. In this case, we say that  $\Psi$ is a L-algebroid morphism.\\

\subsection{Classical derivations on an AL-algebroid}
In  this subsection , $\left(  E,\t,M,\rho,[.,.]_\r \right)  $ will be an AL-algebroid or L-algebroid.

\subsubsection{Lie derivative}

Given any section $s\in \G(\t)$,  the {\bf  Lie derivative} with respect to $s$ on $\Lambda \G^*(\t)$, denoted by  $L_s^{\rho}$,  is the  graded endomorphism, with degree $0$,  characterized by the following properties :

\begin{enumerate}
\item For any function $f\in \Lambda^{0}\G(\t)=\mathcal{F}$
\begin{equation}
L_s^{\rho}\left(  f\right)  =L_{{\rho}\circ s}\left(  f\right)
=i_{{\rho}\circ s}\left(  df\right)  \newline\tag{L0}
\end{equation}
where  $L_{X}$ denote the usual Lie derivative with respect to the vector field  $X$

\item For any $q$--form $\omega\in\Lambda^{q}\G^*(\t)$ (where $q>0$)
\begin{equation}
\left(  L_s^{\rho}\omega\right)  \left(  s_{1},\dots,s_{q}\right)
=L_s^{\r}\left(  \omega\left(  s_{1},\dots,s_{q}\right)  \right)
-{\displaystyle\sum\limits_{i=1}^{q}} \omega\left(  s_{1},\dots,s_{i-1}%
,\left[  s,s_{i}\right]  _{\r},s_{i+1},\dots,s_{q}\right)  \tag{Lq}%
\end{equation}

\end{enumerate}

On the other hand, we can also define for any function $f\in
{\Lambda^{0}\G^*(\t)}=\mathcal{F}$ the element of
$\Lambda^{1}\G^*(\t)$, denoted $d_{\rho}f,$ by %
\begin{equation}
d_{\rho}f={{\rho}^t}\circ df \tag{d0}
\end{equation}

where ${\rho}^t:T^{\ast}M\rightarrow E^{\ast}$  is the transposed mapping of $ \rho$.

\bigskip
 The Lie derivative with respect to $s$ commute with $d_{\rho}$.
\subsubsection{Almost exterior differential}

The {\bf almost exterior differential} on $\L\G^*(\t)$, again  denoted  $d_{\rho}$, (A-differential for short),
is the  graded endomorphism of  degree $1$  characterized by the following properties:
\begin{enumerate}
\item For any function $f\in\ \L^0\G^*(\t)=\mathcal{F}$,
$d_{\rho}f$ is the element of $\L^1\G^*(\t)$
defined by  $d_{\rho}f={{\rho}^t}\circ df $

\item For any $\omega$ in $\L^q\G^*(\t) $
($q>0$), $d_{\rho}\omega$ is the unique element of $\L^{q+1}\G^*(\t)$ such that, for all $s_{0},\dots,s_{q}\in \G(\t)$,
\begin{align*}
\left(  d_{\rho}\omega\right)  \left(  s_{0},\dots,s_{q}\right)   &
={\displaystyle\sum\limits_{i=0}^{q}}\left(  -1\right)  ^{i}L_{s_{i}}^{\rho
}\left(  \omega\left(  s_{0},\dots,\widehat{s_{i}},\dots,s_{q}\right)
\right)  \\
&  +{\displaystyle\sum\limits_{0\leq i<j\leq q}}\left(  -1\right)
^{i+j}\left(  \omega\left(  \left[  s_{i},s_{j}\right]  _{\rho},s_{0}%
,\dots,\widehat{s_{i}},\dots,\widehat{s_{j}},\dots,s_{q}\right)  \right)
\end{align*}
\end{enumerate}
\bigskip

We then have the following properties which are obvious or which can be proved as in finite dimension:
\begin{enumerate}
\item $d_\r(\eta\wedge\zeta)=d_\r(\eta)\wedge \zeta+(-1)^k\eta\wedge d_\r(\zeta)$ for any $\eta\in \L^k\G^*(\t)$ any $\zeta\in \L^l\G^*(\t)$ and any $k,l$ in $\Z$
\item  For a L-algebroid, we have $d_{\rho}\circ d_{\rho}={d_\r}^2=0$. In this case we say that $d_\r$ is the {\bf exterior differential} of the L-algebroid.\\
\end{enumerate}

As in the context of finite dimension (cf. \cite{Ana}), we can prove:

\begin{prop}\label{caractmorphism}${}$\\
Given two AL-algebroids (resp. L-algebroid) $(E_i,\t_i,M,\r_i,[.,.]_{\r_i})$, $\; i=1,2$, let $d_{\r_i}$ be the associated A-differential (resp. exterior differential). For a bundle morphism $\Psi$ from $(E_1,\t_1,M)$ to $(E_2,\t_2,M)$ (over $Id_M$), we denote by $\Psi^*:\L^p\G^*(\t_2)\ap \L^p\G^*(\t_1)$ the induced morphism on $p$-forms. Then, $\Psi$ is an AL (resp. L)-algebroid  morphism if and only if
$$d_{\r_1}  \Psi^*(\o) = \Psi^*(  d_{\r_2}\o)$$
for any $\o\in \L^k\G^*(\t_1)$ and any integer $k>0$.
\end{prop}

Recall that, the bracket $[d_1,d_2]$ of derivations $d_1$ and $d_2$ of the graded algebra $\L\G^*(\t)$ of degree $k_1$ and $k_2$ respectively is the derivation $d_1\circ d_2-(-1)^{k_1k_2}d_2\circ d_1$ of degree $k_1+k_2$. On the graded algebra $\L\G^*(\t)$ with the A-exterior derivation $d_\r$ we have:

\begin{prop}\label{crochetderiv}${}$:\\
For any $s_1$ and $s_2$ in $\G(\t)$, we have
$i_{[s_1,s_2]_\r}(\s)= [[i_{s_1},d_\r],i_{s_2}](\s)$
for any $\s\in \G(\t_*)$
\end{prop}

\noindent \begin{proof}\so{\it Proof}${}$\\

\noindent On one hand, a direct calculation gives the relation

 $ [[i_{s_1},d_\r],i_{s_2}](\s)=L^\r_{s_1}(\s(s_2))-L^\r_{s_2}(\s(s_1))-d_\r(\s)(s_1,s_2)$

 \noindent On the other hand according to the definition of $d_\r$, we get:

$i_{[s_1,s_2]_\r}(\s)=\s([s_1,s_2]_\r)=L^\r_{s_1}(\s(s_2))-L^\r_{s_2}(\s(s_1))-d_\r(\s)(s_1,s_2)$.\\
\end{proof}

\subsubsection{Almost  Schouten-Nijenhuis bracket}
An {\bf almost Schouten-Nijenhuis bracket} (ASN-bracket for short) is  an inner composition law in $\L\G(\t)$
(again) denoted by  $\left[  .,.\right]  _{\rho}$, characterized by the following properties:

\begin{enumerate}
\item $\left[  .,.\right]  _{\rho}$ is a bi-derivation of degree
$-1,$ i.e. an $\R-$bilinear map such that
$\deg\left[  P,Q\right]  _{\rho}=\deg P+\deg Q-1$ which fulfills the following property

$$\left[  P,Q\wedge R\right]  _{\rho}=\left[  P,Q\right]  _{\rho}\wedge
R+\left(  -1\right)  ^{\left(  \deg P+1\right)  \deg Q}Q\wedge\left[
P,R\right]  _{\rho}$$

\item For all $f,g\in\L^0\G(\t)=\cal F$  ,
$\left[  f,g\right]  _{\rho}=0$
\item For all $s\in \L^1\G(\t)=\G(\t)$, $p\in\mathbb{Z}$ and $Q\in\L^p\G(\t)$, $\left[  s,Q\right]  _{\rho}=L_s^{\rho}Q$
\item For all  $s_{1},s_{2}\in\L^1\G(\t)=\G(\t)$, $\left[  s_{1},s_{2}\right]  _{\rho}$ corresponds to the bracket defined on the (A)L-algebroid
\item For all $p,q\in\mathbb{Z}$, $P\in\L^p\G(\t)$, $Q\in\L^q\G(\t)$ ,
$\left[  P,Q\right]  _{\rho}=\left(  -1\right)  ^{pq}\left[  Q,P\right]
_{\rho}$
\end{enumerate}

The ASN-bracket $[.,.]_\r$ is called {\bf Schouten-Nijenhuis bracket} (SN-bracket for short) if, for all $p,q,r\in\mathbb{Z}$ and  $P\in\L^p\G(\t)$,  $Q\in\L^q\G(\t)$\, $R\in\L^r\G(\t)$, the ASN-bracket $[.,.]_\r$ satisfies the graded Jacobi identity:
$$\left(  -1\right)  ^{pr}\left[  \left[  P,Q\right]  _{\rho},R\right]  _{\rho
}+\left(  -1\right)  ^{qp}\left[  \left[  Q,R\right]  _{\rho},P\right]
_{\rho}+\left(  -1\right)  ^{rq}\left[  \left[  R,P\right]  _{\rho},Q\right]
_{\rho}=0$$

Notice that if we take the canonical  L-algebroid $(TM,p_{M},M,Id,[.,.])$  the associated ASN-bracket $[.,.]_{Id}$ is the usual Schouten-Nijenhuis bracket on the graded algebra $\L \G(M)$.

\subsection{Locality of an almost Lie bracket}\label{loc}

In finite dimension it is classical that a
an AL-bracket  $[.,.]_\r$  on an anchored bundle $(E,\t,M\r)$ respects the sheaf of sections of $\t:E\ap M$ or, for short,  is {\bf localizable}  (see for instance \cite{Marl}),  if  the following properties are satisfied:

\begin{enumerate}
\item[(i)]  for any open set $U$ of $M$, there exists a unique  bracket $[.,.]_U$ on the space of sections  $\G(\t_U)$ such that, for any $s_1$ and $s_2$ in  $\G(\t_U)$, we have:
$$[{s_1}_{|U},{s_1}_{|U}]_U=([s_1,s_2]_\r)_{| U}$$
\item[(ii)]  (compatibility with restriction) if $V\subset U$ are open sets, then, $[.,.]_U$  induces a unique AL- bracket  $[.,.]_{UV}$ on $\G(\t_V)$ which coincides with  $[.,.]_V$ (induced by $[.,.]_\r$).
\end{enumerate}

By the same arguments as in finite dimension, when $M$ is smooth regular, we also have:

\begin{prop}\label{localization}${}$:\\
If $M$ is smooth regular then  any  AL-bracket  $[.,.]_\r$  on an anchored bundle $(E,\t,M,\r)$ is localizable.
\end{prop}

If $M$ is not smooth regular, we can no more used the arguments used in the proof of Proposition \ref{localization} . Unfortunately, we have {\bf no example of Lie algebroid} for which the Lie bracket is not localizable. Note that, according to \cite{KrMi}  sections 32.1, 32.4,  33.2 and 35.1,  this problem is similar  to the problem  of localization (in an obvious sense) of  global derivations of  the module of smooth functions on $M$  or  the module  of differential forms on $M$. In \cite{KrMi} and, to our known,  more generally in the literature, there exists no example of such derivations which are not localizable. On the other hand, even if $M$ is not regular, the classical Lie bracket of vector fields on $M$ is localizable.  So, there always exists an anchored bundle ${\cal A}=TM$ and a Lie bracket algebroid $(TM, Id,M,[;,.])$ for which its Lie bracket is localizable.  Moreover,  in Examples \ref{ex1}  we  do not assume   that $M$ is regular but, nevertheless, these Lie brackets are also localizable

\begin{conv}\label{covloc}: in all this paper, from now,  we will assume that either $M$ is smooth regular and if $M$ is not regular, then the Lie bracket $[\;,\;]_\r$ is localizable
\end{conv}

\begin{rem}\label{regparacompact}${}$:
\begin{enumerate}
\item If $M$ is  smoothly paracompact, then $M$ is smooth regular and so any AL-bracket $[.,.]_\r$  on an anchored bundle $(E,\t,M,\r)$ is localizable. On the converse, when $M$ is paracompact,  we can define some AL-bracket $[.,.]_\r$  "locally": given a locally finite covering $\{U_i,\; i\in I\}$ of $M$ and a smooth partition of unity $\{\theta_i,\; i \in I\}$ subordinated to this covering, if $[.,.]_i$ is any AL-bracket on $E_{U_i}$ then:
$$[.,.]_\r=\dis\sum_{i\in I}\theta_i[.,.]_i$$
is an AL-bracket on $(E,\t,M,\r)$. 
\item If $[.,.]_\r$ satisfies the Jacobi identity, then  for any open set $U$,  $[.,.]_U$ satisfies also the Jacobi identity. In this case, $(E_U, \rho_{| U}, U,[.,.]_U)$ is a L-algebroid.
\end{enumerate}
\end{rem}

\begin{rem}\label{intLBalgebroid}${}$\\
Let $\left(  E,\t,M,\rho,[  .,.]_\r  \right)$ be a L-algebroid. As its  bracket is localizable, then ${\cal D}=\r(E)$ is a weak distribution on $M$ (i.e. on each "fiber" ${\cal D}_x$, $x\in M$, we have a Banach structure such that the identity map  from ${\cal D}_x$ as Banach space into ${\cal D}_x$ as normed subspace of $T_xM$ is continuous see \cite{Pel}) . If  the kernel of $\r$ is complemented  in each fiber, then $\cal D$ is integrable (see \cite{Pel}). This situation occurs when $\ker \r$ is finite dimensional or finite codimensional, or when $M$ is a Hilbert manifold.
\end{rem}

\noindent \begin{proof}\so{\it Proof of Proposition \ref{localization} }${}$\\

Let  $s_1 : M\ap  E$ be a smooth section of $\t:E\ap M$ which   vanishes  an open subset $U$ of $M$. We first show that , for any other smooth  section $s_2$  of $ \t: E\ap M$, the bracket
$[s_1,s_2]_\r$ vanishes on $U$. Indeed, choose any point $x\in U$ and  choose a smooth bump function $ f : M \ap \R$  whose support is contained in $U$ and  such that $f (x) = 1$. The section $f s_1$ is a (global section) which  vanishes identically. Therefore, for any other smooth section $s_2$  we have:
$$ 0 = [fs_1,s_2]_\r = -[s_2, fs_1]_\r = -f[s_2,s_1]_\r- df(\r\circ s_2))s_1.$$ 
So at $x$ we have
$$f(x)[s_1,s_2]_\r(x)=df(\r\circ s_2))s_1=0.$$
 Since $f(x) = 1$, we obtain $[s_1,s_2]_\r(x) = 0.$
Now given any open set $U$ in $M$, we must show that the bracket $[\;,\;]_\r$  induces an unique bracket $[\;,\;]_U$ on $\G(\t_{| U})$ such that 
$$[{s_1}_{|U},{s_1}_{|U}]_U=([s_1,s_2]_\r)_{| U}$$
Choose any $x$ in $U$ and, as before, some bump function $ f : M \ap \R$  whose support is contained in $U$ and  such that $f (x) = 1$. Then for any section $s_1$ and $s_2$ in $\G(\t_{| U})$, $fs_1$ and $fs_2$ are global sections of $\t:E\ap M$. So $[fs_1,fs_2]_\r(x)$ is well defined and from the previous argument, this value do not depends of the choice of the bump function $f$. Therefore, we can set
$$[s_1,s_2]_U(x)=[fs_1,fs_2]_\r(x)$$
for any choice of bump function as previously. It follows clearly from this construction that $[\;,\;]_\r$ is localizable.\\
\end{proof}

\noindent \so{\it Local expressions}:${}$\\
Let $(E,\t,M,\r,[.,.]_{\r})$ be an AL-algebroid. In the context of local trivializations (subsection \ref{loctriv}) there exists a field
\begin{center}$
\begin{array}{cccc}
 C: & U & \ap & L(\E\times\E,\E) \\
     & z & \mapsto & C_z
\end{array}
$
\end{center}
\noindent such that 
for $s(z)=(z,u(z))$ and $s'(z)=(z,u'(z))$ we have:
\begin{eqnarray}\label{loctrivct}
[s,s']_U(z)=(z,C_z(u(z),u'(z)))
\end{eqnarray}

 Suppose that the  typical fiber $\E$ has an unconditional basis. According to subsection \ref{locbas},  then for any $x\in M$, there exists an open neighborhood $U$ of $x$ and a set of smooth functions $\{C_{\a\b}^\g,\; \a,\b,\g\in A\}$ on $U$ such that $[.,.]_U$ is characterized by:
\begin{eqnarray}\label{cteloc}
[e_\a,e_\b]_U=\dis\sum_{\g\in A}C_{\a\b}^\g e_\g
\end{eqnarray}
More precisely, if $s=\dis\sum_{\a\in A}s_\a e_\a$ and $r=\dis\sum_{\a\in A}r_\a e_\a$, we have:
\begin{eqnarray}\label{localbracket}
[s,r]_U=\dis\sum_{\a,\b,\g\in A}C_{\a\b}^\g s_\a r_\b e_\g +\sum_{\a\in A}(dr_\a(\r(s))e_\a-ds_\a(\r(r))e_\a)
\end{eqnarray}

\noindent Note that,   the almost exterior differential $d_\rho$ on $\L\G^*(\t_U)$ is also  localizable  i.e. for any open set $U$ in $M$, there exists a unique graded derivation $d_U$ of degree $1$ on  $\L\G^*(\t_U)$, such that $(d_\r\omega)_{| U}=d_U(\omega_{|U})$ and  which is compatible with restriction to open subsets $V\subset U$. 
 So, as for an AL bracket,  in the context of local trivializations (subsection \ref{loctriv}),  given local sections $s(z)=(z,u(z))$ and $s'(z)=(z,u'(z))$ of $E_U$,  and $\o(z)=(z,\xi(z))$  any section of $E^*_U$, according to (\ref{loctrivr}) and  (\ref{loctrivct}) we have
\begin{eqnarray}\label{ds}
d_\r\o(s,s')=<D\xi(\r(s)),u'>-<D\xi(\r(s')),u>- <\xi,[s,s']_\r>\nonumber\\
\textrm{\;\;\;\;\;\;\;\;\;\;\;\;\;}=<D\xi(R(u)),u'>-<D\xi(R(u')),u>-<\xi,C(u,u')>
\end{eqnarray}
where $D\xi$ denotes the differential of the map $\xi:U\ap \E^*$.\\
 
\noindent  In the same way,  the Lie derivative $L_s^\rho$ is  localizable.\\
 
  {\it  For the sake of simplicity, for any open subset $U$ in $M$, we note $[.,.]_\r$,  $d_\r$ and  $L^\r$   instead of its restriction  $([.,.]_r)_{\ U}$, $(d_\r)_U)$ and  $L^\rho_U$, to $U$, respectively}

\section{Almost Lie algebroid defined by a sub AP-morphism or an A-derivation}\label{APmorp}

${}$\indent The notion of sub almost Poisson morphism is a generalization, in the context of Banach manifolds, of  Poisson morphism  or equivalently Lie Poisson structure (see \cite{OdRa2}). The canonical situation corresponds to the weak symplectic structure on the cotangent bundle $T^*M$ of a Banach manifold $M$ (see Example \ref{ex1} {\it 4})

\subsection{Sub almost Poisson morphism }\label{APmor}
Let  $M$ be  a Banach manifold. An {\bf  almost Lie bracket } (AL-bracket for short) on $\cal F$ is a $\R$-bilinear skew-symmetric pairing $\{.,\}$ on $\cal F$ which satisfies the Leibniz rule,  i.e. for any $f,g,h\in {\cal F}$
$$\{f,gh\}=g\{f,h\}+h\{f,g\}.$$

 An {\bf almost Poisson morphism} (AP-morphism for short) on $M$ is a bundle morphism\\ $P:T^*M\ap TM$  which is skew symmetric according to the duality pairing (i.e. such that $<\eta,P\zeta>=-<\zeta,P\eta>$ for any $\eta,\zeta\in T^*M$).

We can associate to such a morphism a $\R$-bilinear skew symmetric pairing $\{.,.\}_P$ on $\G^*(M)$ defined by:
$$\{\eta,\zeta\}_P=<\zeta,P\eta>$$
Moreover, for any $f \in {\cal F}$ we have:
$$\{f\eta,\zeta\}_P= <\zeta,fP\eta>=f\{\eta,\zeta\}_P$$
So,  we get on $\cal F$ an almost Lie bracket $\{.,.\}_P$ defined by
$$\{f,g\}_P=\{df,dg\}_P$$
For any $f\in {\cal F}$  we can associate a unique vector field $grad^P(f)=-P(df)$ which is called the {\bf  almost Hamiltonian vector field} of $f$ (A-Hamiltonian gradient  for short)


Classically, to  an AP-morphism P,  we can associate a skew-symmetric tensor of type $(3,0)$
  $[P,P]  :\G^*(\t)\times\G^*(\t)
\rightarrow\G(\t)$ defined by:

\begin{eqnarray}\label{PP}
\left[  P,P\right] \left(  \eta,\zeta\right)  =P\left(
L_{P\eta}\zeta-L_{P\zeta}\eta+d\left\langle
\eta,P\zeta\right\rangle \right)  +\left[  P\eta,P\zeta\right]
\end{eqnarray}
for all  $\eta,\zeta\in \G^*(\t)$.\\

As in finite dimension, $\{.,.\}_P$ satisfies the  Jacobi identity if and only if the  the tensor $[P,P]$  vanishes identically (see for instance  \cite{MaMo}). In this case, $\cal F$ has a structure of Lie algebra and $(M,\{.,.\}_P)$ is called a {\bf Banach Lie Poisson manifold} (P-manifold for short) (see for instance \cite{OdRa1} or  \cite{OdRa2}). If the Jacobi identity is not satisfied, we say that we have an {\bf  almost Banach Lie Poisson manifold} (an AP-manifold for short). In this case the vector field $grad^P(f)$ is called the hamiltonian gradient of $f$ and  $P$ induces a morphism of Lie algebra between $({\cal F},\{.,.\}_P)$ and $(\G(M),[.,.])$\\

\begin{rem}\label{pbpiosson}${}$\\
 In finite dimension, a Poisson manifold is characterized by a bi-vector $\L$ on a manifold $M$  such that the Schouten-Nijenhuis bracket $[\L,\L]$ vanishes identically.
 On a Banach manifold $M$, an AL-bracket on $\cal F$ gives rise to an element ${\L}$ of $\L T^{**}M$ such that
$${\L}(df,dg)=\{f,g\}$$
Such a bi-vector gives rise to a unique  morphism $P:T^*M\ap T^{**}M$ defined  by the relation:
$$<\zeta,P\eta>=\L(\eta,\zeta).$$
 for all $\eta,\zeta\in \G^*(M)$ where $T^{**}M$ is  the bidual tangent bundle of $M$.\\
 To get a Poisson manifold, we need the additional condition: $P(T^*M) \subset TM$  (see \cite{OdRa1}, \cite{OdRa2}).\\
 Of course  $[{\L},{\L}]=0$ if and only if $[P,P]=0$.
\end{rem}

 For examples and more details about P-manifolds the reader can have a look at  \cite{OdRa1} or  \cite{OdRa2} and some references within these papers.\\

\bigskip

Let $q^{\f}_M:T^{\flat}M\ap M$ be a Banach subbundle of $q_M:T^*M\ap M$. A bundle morphism $P:T^{\flat}M\ap TM$ will be called a {\bf sub almost Poisson morphism}
\footnote{ this terminology is chosen in analogy to sub-riemannian structures  on a manifold }
 (sub AP-morphism for short) if $P$ is skew-symmetric  relatively to the duality pairing i.e.  $<\a,P\b>=-<\b,P\a>$  for any $\a, \b \in T^\f M$. As before, we get  a $\R$-bilinear skew-symmetric pairing $\{.,.\}$ on $\G(q^\f_M)$. The set  ${\cal F}^\f=\{f\in {\cal F} :  df\in \G(q^\f_M)\}$ is a sub-algebra  of $\cal F$. Using the same arguments as in \cite{KrMi}, section 48, we can see that,  as above, $P$ induces on  ${\cal F}^{\flat}$ an almost Lie bracket which will be again denoted by $\{.,.\}_P$. We will say that $(M, {\cal F}^{\flat},\{.,.\}_P)$ is a {\bf  sub almost Banach Lie Poisson} manifold (sub AP-manifold for short). Of course, when $\{.,.\}_P$ satisfies the Jacobi identity, $({\cal F}^{\flat},\{.,.\}_P)$,  has a Lie algebra structure and we say that ${({\cal F}^{\flat},\{.,.\}_P})$ is a {\bf sub Banach Lie-Poisson} manifold  (sub P-manifold for short). Then for any $f\in{\cal F}^\f$ we can associate a {\bf sub  almost Hamiltonian gradient} $grad^P(f)=-P(df)$ (sub A-Hamiltonian gradient for short) \\

Of course, to a sub AP-morphism $P$ on $M$, we can also  associate a skew symmetric tensor of type $(3,0)$: $[P,P]:\G(q^\f_M)\times\G(q^\f_M)
\rightarrow\G(\t)$ using the same definition as in (\ref{PP}), but for $\a,\b\in T^{\flat}M$.  Again, $\{.,.\}_P$ satisfies the Jacobi identity if and only if $[P,P]$ vanishes identically. In this case, $({\cal F}^\f,\{.,.\}_P)$ has a Lie algebra structure and $P$ induces a Lie algebra morphism from $({\cal F}^\f,\{.,.\}_P)$ to $(\G(M),[.,.])$\\

Consider two sub AP-morphisms $P_i$ on the manifolds $M_i$, $\;i=1,2$ and   a map $\psi:M_1\ap M_2$. We say that $(M_1,P_1)$ and $(M_2,P_2)$ are ${\bf \psi}${\bf -related} if we have:
$$P_2= (T\psi)^*\circ P_1\circ T\psi$$
In this case for the associated AL-bracket $\{.,.\}_{P_i}$ on the algebra ${\cal F}^\f_i$ $i=1,2$ we have:
$$\{f\circ\psi,g\circ \psi\}_2=\psi\circ\{f,g\}_1$$
for all $f,g\in {\cal F}^\f_1$. Moreover, in this case, $P_1$ is a sub P-morphism if and only if $P_2$ is a sub P-morphism  and then $\psi$ gives rise to a Lie algebra morphism between $({\cal F}^\f_1,\{.,.\}_1)$ and $({\cal F}^\f_2,\{.,.\}_2)$.\\

A lot of results about AP-manifolds and P-manifolds can be extended to the context of sub AP-manifolds and sub-P-manifolds when considering the structure of Lie algebra of $({\cal F}^{\flat},\{.,.\}_P)$. We do not develop these aspects here.

\begin{exe}\label{ex2}${}$\\
 Let $\o$  be a non degenerated $2$-form on a manifold $M$. We denote by $\o^\f: TM\ap T^*M$ the associated morphism defined by $\o^\f(X)=i_X \o$ for $X\in T_xM$. Suppose that $T^\f M=\o^\f(TM)$ is a Banach subbundle of $T^*M$. Then $\o^\f$ is an isomorphism from $TM$ onto $T^\f M$. So $P= (\o^\f)^{-1}: T^\f M\ap TM$
 is a sub AP-morphism. We get  a sub P-morphism if and only if $\o$ is closed.
 In particular, if $M$ is the cotangent bundle $T^*N$ of a Banach manifold $N$, if $\o$ is the $2$ fundamental form on $T^*N$, we get a sub Poisson structure  on $T^*N$ which corresponds to the natural weak symplectic structure on $T^*N$ (see Example \ref{ex1} {\it 4}). Note that we get a Poisson structure on $T^*N$ if and only if  $ N$ is modeled on a reflexive Banach space. \\
 \end{exe}

\begin{rem}${}$\\
  Consider  a sub-Poisson morphism $P:T^{\flat}M\ap TM$, and denote by ${\cal D} = P(T^{\flat}M)$ the associated (weak) distribution on $M$. If the kernel of $P$ is complemented in each fiber, then $\cal D$ is integrable and each leaf is a weak symplectic manifold (see \cite{Pel}). This situation is always satisfied when $\ker P$ is finite dimensional or finite co-dimensional (for instance if $P$ is Fredholm, or injective) , or when $M$ is an Hilbert manifold.\\
\end{rem}

\bigskip
Let $\t :E\ap M$ be a Banach bundle and $\t_*:E^*\ap M$ its associated dual Banach bundle. We denote by ${\cal F}(E^*)$ the set of smooth functions on $E^*$. Any $f\in {\cal F}(E^*)$ is called {\bf  linear}, if the restriction of $f$ to each fiber $E_x^*=\t_{*}^{-1}(x)$ is a linear map.\\


\begin{defi}${}$\\
 Let $q^{\f}:T^{\flat}E^*\ap E^*$ be a subbundle of  $q_{E^*}:T^*E^*\ap E^*$. Consider a sub AP-morphism\\ $P:T^{\flat}E^*\ap TE^*$ and its associated bracket $\{.,.\}_P$ on ${\cal F}^{\flat}(E^*)\subset  {\cal F}(E^*)$.
 The sub AL-morphism $P$ on $E^*$  is called linear if for any linear functions $f$ and $g$ on $E^*$  which belong to ${\cal F}^{\flat}(E^*)$ their bracket $\{f,g\}_p$ is linear.
\end{defi}

\begin{exe}\label{ex3}${}$\\
In the context of Example \ref{ex2},  if we take  $E=TM$, the sub P-morphism\\ $\Pi=(\o^\f)^{-1}:T^\f (T^*M)\ap T(T^*M)$ is a linear sub P-morphism on $T^*M$.\\
\end{exe}

${}$

 To a given  Banach bundle $\t :E\ap M$, {\it  we will construct a canonical subbundle $T^{\f}E^*$ on which will be defined all "interesting" sub AP-morphisms} (see Proposition \ref{admis} and subsection \ref{alg-Poisson}).\\

${}$
First of all,   for any section $s\in \G(\t)$, we  associate the linear function $\Phi_s$ on $E^*$ defined by
$$\Phi_s(\xi)=<\xi,s\circ\t_*(\xi)>$$
We then have  the following properties:

 \begin{lem}\label{ecritlocale}${}$
 \begin{enumerate}
 \item The map $s \mapsto \Phi_s$ is linear and injective; we also have $\Phi_{fs}=(f\circ \t_*)\Phi_s$ for any $s\in \G(\t)$ and $f\in {\cal F}$
 \item In  local trivializations (subsection \ref{loctriv}) we have

 $<d\Phi_s(\s),w_2> =<w_2,s\circ\t_*(\s)>$ for any $ w_2\in \E^*$, considered  as vertical fiber $T_\s E^*$;

  if  for some $f\in {\cal F}$, we have $\; d(\Phi_s+f\circ \t_*)(\s)=0$ then $s\circ \t_*(\s)=0$.\\
\end{enumerate}
\end{lem}

\noindent \begin{proof}\so{\it Proof of Lemma \ref{ecritlocale}}\\

The first part is easy and left to the reader.

With the previous notations, if $w_2$ belongs to the vertical part of $T_\s
E^*$\ which is  $\{\s\}\times E^*$ with our notations. We can consider $w_2$ as an element of $E^*_{\t_*(\s)}$. We then have

$<d\Phi_s(\s),w_2>=\dis\lim_{t\ap 0}1/t[<\s+tw_2,s\circ\t_*(\s)>-<\s,s\circ\t_*(\s)>]=<w_2,s\circ\t_*(\s)>$.

Assume that $d(\Phi_s+f\circ \t_*)(\s)=0$. Given any $w_2\in  E^*_{\t_*(\s)}$. As before, $w_2$ can be considered as a vector in the vertical part of $T_\s E^*$. From the previous relation we get:

$<d\Phi_s(\s),w_2>=-<df(\t_*(\s)), T_\s\t_*(w_2)>=0$

which ends the proof.

\end{proof}

\begin{prop}\label{admis}${}$\\
The set
 $$T^{\f}E^*=\dis{\bigcup}_{x\in M}\{(\s,\eta)\in T_\s^*E^*,  \;\; \eta=d(\Phi_s+f\circ\t_*),\; \; s\in \G(\t_U),\;\; f\in {\cal F}(U), \;\; U \textrm{ neighbourhood of } x \}$$
 is a well defined subset of $T^*E^*$ and if $q^{\f}$ is  the restriction of  $q_{E^*}$ to $T^{\flat}E^*$ then $$q^{\f}: T^{\flat}E^*\ap M$$ is a  Banach bundle of typical fiber $\M^*\times\E$. In particular, $T^{\flat}E^*=T^*E^*$ if and only if $\E$ is reflexive.
\end{prop}

\noindent \begin{proof}\so{\it Proof of Proposition \ref{admis}}\\

First notice that, using our notations, at any $\s\in E^*$, any vector $(\s,\eta_1)\in \{\s\}\times\M^*$ can be written as $\eta_1=df(\t_*(\s))$ for some $f\in{\cal F}$ and moreover, this choice only depends on the value $df(\t_*(\s))$. On the other hand, from Lemma \ref{ecritlocale} part {\it 2}, we can identified the restriction of $d\Phi_s(\s)$ to the vertical part to $T_\s E^*$ with $s\circ\t_*$.
 So,  the set
 $$\{d\Phi_s(\s)_{|\{\s\}\times \E^*},\;\; s\in \G(\t_U)\}$$  generates the subspace $\{\s\}\times\E$, considered as subspace of the vertical part $\{\s\}\times  \E^{**}$ of $T_\s^*E^*$.\\

Assume that we can write $\eta=d(\Phi_{s_1}+f_1\circ\t_*)(\s)=d(\Phi_{s_2}+f_2\circ\t_*)(\s)$.
From Lemma \ref{ecritlocale} part {\it 2}, we get $s_1\circ\t_*(\s)=s_2\circ\t_*(\s)$ and  then we must have $df_1(\t_*(\s))=df_2(\t_*(\s))$. It follows that the subset
$$T_\s^{\flat}E^*=\{(\s,\eta)\in T_\s^*E^*,  \;\; \eta=d(\Phi_s+f\circ\t_*),\; \; s\in \G(\t_U),\;\; f\in {\cal F}(U), \;\;U \textrm{ neighbourhood of } x\}$$
is well defined and is a subspace of $T_\s^*E^*$  which, with our notations, is exactly $\{\s\}\times \M^*\times \E$. From the local definition of   $T_\s^{\flat}E^*$, it follows that the restriction  $q^{\flat}$ of  $q_{E^*}: T^{*}E^*\ap E^*$  gives rise to a Banach subbundle. Of course, we get $T^{\f}E^*=T^*E^*$ if  and only if $\E$ is reflexive.\\
\end{proof}

\subsection{Relation between AP-algebroid and sub AP-morphism }\label{alg-Poisson}

Now we are able to give an adaptation to the Banach context of the classical result about equivalence between Poisson structure  on $E^*$ and structure of AL-algebroid on $E$ (see for example  \cite{Marl}).

\begin{theo}\label{linearAP}${}$\\
Let  $P: T^{\flat}E^*\ap TE^*$ be a linear sub AP-morphism  on $E^*$. Then there exists a unique   AL-algebroid structure  $(E,\t,M,\r,[.,.]_P)$ characterized by:
\begin{eqnarray}
\label{linAP1} \Phi_{[s_1,s_2]_P}=\{\Phi_{s_1},\Phi_{s_2}\}_P,\; \textrm{ for any } s_1,s_2\in \G(\t)\\
\label{linAP2} \{\Phi_s,f\circ \t_*\}_P=L^\r_s(f)\circ \t_*,\; \textrm{ for any } f\in {\cal F},\; s\in \G(\t)
\end{eqnarray}
Moreover, $(E,\t,M,\r,[.,.]_P)$ is a   L-algebroid if and only if $P$ is a sub-Poisson morphism.\\
Conversely, each  AL-algebroid structure $(E,\t,M,\r,[.,.]_\r)$ defines a unique linear  bracket $[.,.]_\r$ on the sub-ring  ${\cal F}^{\flat}(E^*)$ which is associated to a unique linear sub AP-morphism   on $ E^*$ which is characterized by relations  (\ref{linAP1}) and (\ref{linAP2}).  Moreover,  $P$ is a sub-Poisson morphism on $E^*$ if and only if $(E,\t,M,\r,[.,.]_\r)$ has a L-algebroid structure.
\end{theo}

\begin{exe}\label{ex4}${}$\\
According to Example \ref{ex2} and Example \ref{ex3}, as in finite dimension, the canonical sub P-morphism $\Pi$ induces on the  bundle $(TM,p_M,M)$ a structure of L-algebroid which is, in fact, the canonical  L-algebroid structure (see Example \ref{ex1} {\it 1}).
\end{exe}
\bigskip
\noindent \so{\it Local expressions} : \\

 In the context of local trivialization (subsection \ref{loctriv}), recall that $TE^*_U\equiv U\times \E^*\times\M\times \E^*$ and $T^\f E^*_U\equiv U\times\E^*\times\M^*\times\E$. Given an AP-morphism $P$, recall  that  locally we have the following characterization (see (\ref{loctrivr}) and (\ref{loctrivct}))

 $\r(s)=(., R(u))$ for any $s(z)=(z,u(z))$

 $[s,s']_\r=(\;,C(u,u'),)$ for any $s(z)=(z,u(z))$ and $s'(z)=(z,u'(z))$

 \noindent  so locally we have
   \begin{eqnarray}\label{Ploc}
   \begin{array}{ll}
     <(df',u'),P(\o(x))(df,u)> & =<\xi,[s,s']_\r>(x)+ <df', \r(x,u)>- <df, \r(x,u')> \\
      &=<\xi,C(u,u')>+<df',R(u)>-<df,R(u')>
   \end{array}
\end{eqnarray}
for any functions $f$ and $f'$ on $U$, any sections $s$ and $s'$ any form $\o(z)=(z,\xi(z))$.\\

 If $F$ is a  function in  ${\cal F}^\f(E^*_U)$, in local trivialization, we can write $dF=(D_1F,D_2F)$ on $U\times \E^*$ where $D_1F$ (resp. $D_2F$) is the partial derivative of $F$ according to the first  factor (resp. the second factor). Notice that, $D_2F(\s)$ belongs to $\E$. So for any $F,G\in {\cal F}^\f(E^*_U)$, their  Poisson bracket is given by:
 \begin{eqnarray}\label{poisbracket}
\{F,G\}_P(\s)= <D_1F,R(D_2G)>(\s)\;-\;<D_1G,R(D_2F)>(\s)\;-\;<\xi,[D_2F,D_2G]_\r>(\s)
\end{eqnarray}

\noindent Suppose that each  Banach  space $\E$  and $\M$ has an unconditional basis. According to subsection \ref{locbas} and also   local trivializations (subsection \ref{loctriv}),  recall that  we have the following local coordinates:

$\bullet\;\;$ $(x,u)=(x^i,u^\a)$ on $E_U$

$\bullet\;\;$   $\s=(x,\xi)=(x^i,\xi_\a)$  on $E^*_U$ ("weak-* coordinates" for $\xi_\a$ )

$\bullet\;\;$ on $E^*_U$,  the  tangent space  $T_\s E^*_U$ is generated by the basis $\{\dis\frac{\partial}{\partial x^i} \}_{i\in I}$ and "weakly-*"

generated by the basis  $\{ \dis \frac{\partial}{\partial \xi_\a}\}_{\a\in A}$;

$\bullet\;\;$ on $E^*_U$, according to Lemma \ref{ecritlocale} part 1 and its proof, each fiber $T_\s^{\f}E^*$ is generated by the basis $\{dx^i\}_{i\in I}$ and $\{e_\a\circ \t_*\}_{\a\in A}$.
Notice that, we have $d\xi_\a=e_\a\circ \t_*$.\\

\noindent With  these notations, any sub AP-morphism $P:T^{\f}E^*\ap TE^*$ associated to an AL-bracket $[.,.]_\r$ (as in Theorem \ref{linearAP}) are related  in the following way:

the Lie bracket $[.,.]_P$ is locally characterized by  $[e_\a,e_\b]_U=\dis\sum_{\g\in A}C_{\a\b}^\g e_\g$ (see (\ref{cteloc}));

the anchor $\r$ can be characterized by $\r(e_\a)=\dis\sum_{i\in I}\r^\a_i\frac{\partial}{\partial x^i}$ (see (\ref{locrho}))

the AP-morphism $P$ is characterized by:
\newline

$P(d x^i)=\dis\sum_{\a\in A}\r_\a ^i \frac{\partial}{\partial \xi_\a} \hfill{}$

$P(d\xi_\a)=-\dis\sum_{i\in I}\r_\a ^i \frac{\partial}{\partial x^i}-\dis\sum_{\b,\g\in A} C_{\a \b}^\g \xi_\g \frac{\partial}{\partial \xi_\b}$\\

\noindent On $E^*_U$,  the associated algebra of functions  ${\cal F}^\f(E^*_U)$ are functions on $E^*_U$ which depend on variables $(x_i)_{i\in I}$ and $(\xi_\a)_{\a\in A}$. The AL-bracket $\{.,.\}_P$ on ${\cal F}^\f(E^*_U)$ is characterized by:

$\{\xi_\a,\xi_\b\}_P=-\dis\sum_{\g\in A}C_{\a\b}^\g\xi_\g,\;\;\; \{x^i,\xi_\a\}_P=\r^i_\a,\;\;\;\; \{x^i,x^j\}_P=0$

The Poisson Bracket of  $F,G\in {\cal F}^\f(E^*_U)$ is given by:

$\{F,G\}_P=\dis\sum_{i\in I,\a\in A}\r^i_\a(\dis\frac{\partial F}{\partial x^i}\dis\frac{\partial G}{\partial \xi_\a}-\dis\frac{\partial G}{\partial x^i}\dis\frac{\partial F}{\partial \xi_\a})-\dis\sum_{\a,\b,\g\in A}C_{\a\b}^\g\xi_\g\dis\frac{\partial F}{\partial \xi_\a}\dis\frac{\partial G}{\partial \xi_\b}$

\bigskip\bigskip

\noindent For the first part of the proof   of Theorem \ref{linearAP}, we need the following Lemma:

\begin{lem}\label{proplinearAP}${}$\\
Consider a linear sub AP-morphism $P: T^{\flat}E^*\ap TE^*$ on $E^*$ and $\{.,.\}_p$ the associated bracket on ${\cal F}^{\flat}(E^*)$.
\begin{enumerate}
 \item  for any section $s\in \G(\t)$ and any $f\in{\cal F}$ the bracket $\{\Phi_s,f \circ\t_* \}_P$ belongs to $\cal F$ .
\item If $f $ and $g$  belong to $\cal F$,  then $\{f\circ \t_*,g\circ \t_*\}_P =0$.
\end{enumerate}
\end{lem}

\noindent\begin{proof}\so{\it Proof of Lemma \ref{proplinearAP}}\\
\indent  This  proof is an adaptation of  the proof of  the analogous result  in finite dimension, (cf  \cite{LMM} for instance)\\

Consider any $s, s' \in \G(\t)$, and $f\in {\cal F}$. Using the Leibniz rule for the bracket $\{.,.\}_P$ we get:
\begin{eqnarray}\label{dvptlin}
\{\Phi_{s},(f\circ \t_*)\Phi_{s'}\}_P =(f\circ \t_*)\{\Phi_{s},\Phi_{s'}\}_P +\Phi_{s'}\{\Phi_s,f\circ\t_* \}_P
\end{eqnarray}
Moreover  $\{\Phi_s,(f\circ \t_*)\Phi_{s'}\}_P$ is a linear function on $E^*$. Since $(f\circ \t_*)\{\Phi_s,\Phi_{s'}\}_P$ is also a linear map, from (\ref{dvptlin}) it follows that $\Phi_{s'}\{\Phi_s,f\circ\t_* \}_P$ is a linear function on $E^*$ for any $s'\in \G(\t)$. So $\{\Phi_s,f\circ\t_* \}_P$ must be constant on each fiber and then   we get {\it 1}.\\
On the other hand, by same argument, we have
$$ \{(f\circ\t_*)\Phi_{s'},g\circ\t_*\}_P =(f\circ\t_*)\{\Phi_{s'},g\circ\t_*\}_P +\Phi_{s'}\{f\circ\t_*,g\circ\t_*\}_P$$
So  from  part {\it 1}, we deduce that , $\Phi_{s'}\{f\circ\t_*,g\circ\t_*\}_P$ belongs to $\cal F$ for any $s'\in \G(\t)$. If follows that we must have  $\{f\circ\t_*,g\circ\t_*\}_P=0$ and we get {\it 2}.\\
\end{proof}

 \noindent\begin{proof}\so{\it Proof of Theorem \ref{linearAP}} (adaptation, in our context,  of the proof of \cite{LMM} of the same result in finite dimension)\\

 \noindent  Consider a linear sub AP-morphism  $P: T^{\flat}E^*\ap TE^*$ on $E^*$ and $\{.,.\}_P$ the associated bracket on ${\cal F}^{\flat}(E^*)$. As $\Phi$ is injective, the pairing $[.,.]_P$ defined by  (\ref{linAP1}) is  well defined and is $\R$-bilinear and skewsymmetric.  From Lemma \ref{proplinearAP} part {\it 1} and  the Leibniz property of  $\{.,.\}_P$  it follows that, for some fixed $s\in \G(\t)$, the map $\r(s):f \mapsto  \{\Phi_s, f\circ \t_*\}_P$  defines  a derivation on ${\cal F}$. On the other hand, we have:
 $$ \{\Phi_s, f\circ \t_*\}_P=<d(f\circ\t_*),P(d\Phi_s)>=<df,d\t_*\circ P(d\Phi_s)>$$
It follows that $\r(s)=d\t_*\circ P(d\Phi_s)$ is a vector field on $M$. Notice that, from the properties of $\Phi$ (Lemma \ref{ecritlocale} part {\it 1}), the Leibniz  property of  $\{.,.\}_P$, and Lemma \ref{proplinearAP} part {\it 2}, we have:
 \begin{eqnarray}\label{morpr}
{\r}(f s)=f{\r}(s)
\end{eqnarray}
for any $f\in {\cal F}$. As $P$ is a bundle morphism, it follows that the bracket $\{.,.\}_P$   is localizable, so on one hand the same is true for the bracket $[.,.]_P$ and on the other hand, from (\ref{morpr}) we get that ${\r}(s)$ only depends on the value of  $s$ at any point $x\in M$ so we get a morphism bundle $\r:E\ap M$  defined by $\r(u)={\r}(s)(x)$ where $s$ is any (local) section such that $s(x)=u$.

From  (\ref{linAP1}) and (\ref{linAP2}), and Lemma \ref{ecritlocale} part {\it 1}, for any $s_1,s_2\in \G(\t)$ and $f\in {\cal F}$ we have:

$$\Phi_{[s_1,f s_2]}=  \{\Phi_{s_1},(f\circ\t_*)\Phi_{s_2}\}_P=(f\circ\t_*)\{\Phi_{s_1},\Phi_{s_2}\}_P+\Phi_{s_2}L_{\r(s_1)}(f)\circ\t_*$$
So we get the following relation:
$$[s_1,fs_2]_P=f[s_1,s_2]_P +\r(s_1)(f)s_2$$
Thus we have proved that $(E,\t,M,\r,[.,.]_P)$ is an AL-algebroid. From the properties of $\Phi$, it follows that if  $[.,.]_P$ satisfies the Jacobi identity, it implies that $\{.,.\}_P$ condition is true  if and only if $P$ is a Poisson morphism.\\

Conversely,  let $(E,\t,M,\r,[.,.]_\r)$ be an AL-algebroid. Note that, from convention \ref{convloc} the bracket $[\;,\;]_\r$ is localizable. We want to associate a linear  AP-morphism $P:T^*E^*\ap TE^*$ such that, according to the first part,  the induced AL-algebroid structure on $(E,\t,M)$ is exactly $(E,\t,M,\r,[.,.]_\r)$. For this  we must have:

$\{\Phi_s,f\circ\t_*\}_P=L_{\r\circ s}(f)\circ\t_*\textrm{ and } \{f\circ\t_*,g\circ\t_*\}_P=0$

\noindent Locally, with the  notations of Lemma \ref{ecritlocale}, we define $P_\s:T_\s^{\flat}E^*\ap T_\s E^*$ as follows:

for  $\eta=d(\Phi_s+f\circ \t_*)(\s)$,  and $\eta'=d(\Phi_s'+f'\circ \t_*)(\s)$
we set
\begin{eqnarray}\label{L}
\L(\s)(\eta,\eta')=\Phi_{[s,s']_\r}(\s)+L_{\r (s)} (f')\circ\t_*(\s)-L_{\r(s')}(f)\circ\t_*(\s)
\end{eqnarray}
 As we have seen that $\eta$ (resp. $\eta'$) only depends on  $s\circ \t_*(\s)$ and $df\circ \t_*(\s)$ (resp. $s' \circ \t_*(\s)$ and $df' \circ \t_*(\s)$), then it follows that $\L$ is well defined at $\s$. Moreover,  from its local definition, we get a smooth section of $\L^2 T^\f E^*$. It follows that the map $\eta \mapsto \L(\s)(.,\eta)$ defines a linear map  $P(\s):T^\f_\s E^* \ap [T^\f_\s E^*]^*$ directly given by:
  \begin{eqnarray}\label{P}
P_\s(\eta)=\Phi_{[s,.]_\r}(\s)+L_{\r (s)} (.)\circ\t_*(\s)-L_{\r(.)}(f)\circ\t_*(\s)
\end{eqnarray}

\noindent Using our notations, we have $T^\f_\s E^*\equiv\{\s\}\times \M^*\times \E$ and so $[T^\f_\s\E^*]^*\equiv \{\s\}\times\M^{**}\times \E^*$.


\noindent Recall that  $\eta=d\Phi_s +df\circ \t_*$. On one hand, any $\o\in \{\s\}\times \M^*$ can be written as $\o=dg\circ \t_*(\s)$ for some $g:U\subset \M\ap \R$; so we have

$<\o,P_\s(\eta)>=\L_\s(\eta,dg\circ\t_*)=L_{\r (s)} (g)\circ\t_*(\s)$

\noindent It follows that  $P_\s(\eta)_{| \{\s\}\times \M^*}$  belongs to $\{\s\}\times \M$ considered as a  subspace of $\{\s\}\times \M^{**}$.

On the other hand we have
$$
\begin{array}{ll}
  <d\Phi_{s'},P_\s(\eta)>=\L_\s(\eta,d\Phi_{s'}) & =\L_\s(d\Phi_{s},d\Phi_{s'})+\L(\s)(df\circ\t_*,d\Phi_{s'}) \\
   & =\Phi_{[s,s']}(\s)- L_{\r(s')}(f)\circ \t_*(\s)
\end{array}
$$
\noindent It follows that   $P(\s)(\eta)_{| \{\s\}\times \E}$ belongs to $\{\s\}\times \E^*\times \M$. In conclusion, $P(\s)(\eta)$ belongs to\\ $T_\s E^*\equiv \{\s\}\times \M^*\times E$.

Notice that  the definition (\ref{P}) of $P$ is in fact local so $P$ is a smooth bundle morphism. \\

Of course as usually,   $P$  gives rise to a bracket on ${\cal F}^\f(E^*)$ which  is exactly given by (\ref{L}). This bracket is denoted by $\{.,.\}_P$.  As $P$ is a bundle morphism, $\{.,.\}_P$ is localizable. Moreover, it  satisfies the relations
\begin{eqnarray}
\{\Phi_{s_1},\Phi_{s_2}\}_P= \Phi_{[s_1,s_2]_\r}=<d\Phi_{s_2},P(d\Phi_{s_1})>,\; \textrm{ for any } s_1,s_2\in \G(\t)\label{Phi}\\
L_{\r(s)}(f)\circ \t_*= \{\Phi_s,f\circ \t_*\}_P=<df\circ\t_*, P(d\Phi_s)>,\; \textrm{ for any } f\in {\cal F},\; s\in \G(\t)\label{Phif}\\
\{f_1\circ\t_*,f_2\circ\t_*\}_P=<df_2\circ\t_*,P(df_1)\circ \t_*>=0,\; \textrm{ for any } f _1, f_2\in {\cal F}\label{f}.
\end{eqnarray}


Let $F$ be a smooth linear function on $E^*$ which belongs to ${\cal F}^\f$. Fix a point $\s=(x,\xi)\in E$ and $E^*_U\equiv U\times \E^*$ a neighborhood of $\s$. We denote by $d_2 F$ the  partial differential of $F$ relative to $\E^*$ in the product $U\times\E^*$. As $F$ is linear,  there exists a section $S:U\ap \E^{**}$ such that
$$d_2F(\zeta)=<S,\zeta>$$
for any $\zeta \in \E^*$.
But  $F$ belongs to  ${\cal F}^\f$ so $dF$ is a section of $T^\f E^*\ap E^*$; then we must have $S:U\ap E$ and $d_2F=d_2\Phi_S$. Then, from (\ref{Phi}) it follows that $P$ is a linear AP-morphism. On the other hand, the previous relations (\ref{Phi}), (\ref{Phif}) and (\ref{f})  mean that the Lie bracket $[.,.]_P$ induced by $P$ on $\G(\t)$ is exactly the original one  $[.,.]_\r$. Finally, if $[.,.]_\r$ satisfies the Jacobi identity, the previous relation implies that $\{.,.\}_P$ also satisfies the Jacobi identity for functions of type $\Phi_s+f\circ\t_*$. So it follows that $[P,P]$ vanishes indentically and then, the Jacobi identity is satisfied for any $f,g,h\in {\cal F}^{\flat}(E^*)$, which ends the proof.\\
\end{proof}

\subsection{Relation between AP-algebroid  and  A-derivations}\label{A-deriv}

Recall that an  {\bf A-exterior differential} is a graded  derivation $\d$ of degree $1$ of $\L\G^*(\t)$ which is {\bf localizable } i.e. for any open set $U$ in $M$, there exists a unique graded derivation $\d_U$ of degree $1$ of $\L\G^*(\t_U)$ such that
$$(\d\omega)_{| U}=\d_U(\omega_{|U})$$
which is compatible with restriction to open subsets $V\subset U$  and satisfies the following properties:

\begin{enumerate}
\item $\d(\eta\wedge\zeta)=\d(\eta)\wedge \zeta+(-1)^k\eta\wedge \d(\zeta) $ for any $ \eta\in \L^k\G^*(\t) \textrm{ any } \zeta\in \L^l\G^*(\t)\textrm{ and any } k,l \in \Z$
\item For a L-algebroid, we have $\d\circ \d=\d^2=0$. In this case we say that $\d$ is a  {\bf exterior differential}.
\end{enumerate}

We will adapt the classic result obtained in finite dimension: given a derivation $\d$, one  associates a unique bracket $[.,.]_\d$ on $\G(\t)$ with anchor $\r$ such that the almost  exterior derivative $d_\r$ associated is exactly $\d$.
However, in finite dimension,  any (local) derivation of $\cal F$ is a vector field, but in infinite Banach context it is not true (see subsection \ref{deriv}). So, we must impose another condition on $\d$ to get an analogous  result. Moreover, in finite dimension, the exterior algebra $\L\G^*(\t)$ is locally generated by its elements of degrees $0$ and $1$ which is not true in the Banach framework, even when we have Schauder basis (see Remark \ref{pbderiv}).\\

  In the context of the Proposition \ref{crochetderiv}  first we have:

\begin{lem}\label{conddelta}${}$\\
Consider a graded A-derivation $\d$ of degree $1$ of $\L\G^*(\t)$ which is  localizable. For any $s_1$ and $s_2$ in $\G(\t)$,  the bracket $[[i_{s_1},\d],i_{s_2}]$ is a derivation of degree $-1$ and its restriction to $\G(\t_*)$ can be identified with a section of the bi-dual $E^{**}\ap M$
\end{lem}

\noindent\begin{proof}\so{\it Proof}${}$\\

First by construction, the degree of $D=[[i_{s_1},\d],i_{s_2}]$  is $-1$ so $D$ maps $\L^1\G(\t_*) $ into $\cal F$ and $Df=0$ for any $f\in {\cal F}$. It follows that $D(f\s)=fD\s$ for any $\s\in \G(\t_*)$ and $f\in {\cal F}$ and so the map $\s \mapsto D\s $ is ${\cal F}$ is linear, which ends the proof.

\end{proof}

Now, in our context, we have
\begin{theo}\label{algebroidetdiff}${}$\\
Let  $\t :E\rightarrow M$ be a Banach bundle and  $\t_*:E^{\ast}\rightarrow M$  its dual bundle. Consider  a localizable graded derivation $\delta$ of degree $1$ of the graded algebra  $\L\G^*(\t)$. Assume that  for any $s_1$ and $s_2$ in $\G(\t)$,  the bracket $[[i_{s_1},\d],i_{s_2}]$ in restriction to $\G(\t_*)$ can be identified with a section of $E\subset E^{**}\ap M$

Then $\delta$ defines a unique  bracket  $\left[.,.\right]_\d$ on $\G(\t)$ and there exists a unique morphism $\r:E\ap M$ such that:
\begin{enumerate}
  \item  for any function $f\in\mathcal{F}$ and any section $s\in \G(\t)$  we have
  \begin{eqnarray}\label{rhodif}
\rho ( s)(f)=<\d f,s>
\end{eqnarray}
 \item $[.,.]_\d$ is characterized  by
  \begin{eqnarray}\label{bracketdif}
 \s([ s_{1},s_{2}]  _{\d})=\d(\s(s_1))(s_2)-\d(\s(s_2))(s_1)-\d\s(s_1,s_2), \textrm{ for any }  \s \in \G(\t_*)
\end{eqnarray}
\end{enumerate}

\noindent In particular $(E,\t,M,[.,.]_\d)$ is an AL-algebroid and the almost exterior derivative associated to this structure coincides with  $\delta$ on $\L^0\G^*(\t)=\cal F$ and $ \L^1\G^*(\t)=\G(\t_*)$.
Moreover, $(E,\t,M,[.,.]_\d)$ is a L-algebroid if and only if $\d^2=0$
\end{theo}

\noindent \begin{proof}\so{\it Proof}.\\
From (\ref{rhodif}) and our assumption, we get a linear map $\r:\G(\t)\ap \G(M)$ which gives rise to a linear map $\r_U:\G(\t_U)\ap \G(U)$ for any open $U$ in $M$. On the other hand, as  $\d f$ is a $1$-form on $E$, for any smooth function $h$ defined on an open $U$ we have: $<\d f_{| U},hs>=h<\d f_{| U},s>$; it follows that $\r(s)$ only depends on the value of $s$ at each point. So, we get a bundle morphism from $E$ to $TM$.

Notice that, according to the proof of Proposition  \ref{crochetderiv} the RHS of (\ref{bracketdif}) is exactly $[[i_{s_1},\d],i_{s_2}(\s)]$. Taking into account the definition of the map $s \mapsto \Phi_s$ in subsection \ref{APmor}, the LHS of (\ref{bracketdif}) is exactly $\Phi_{[s_1,s_2]_\d}(\s)$. From our assumption, there is a section which we can denote by  $[s_1,s_2]_\d$ of $\t:E\ap M$ such that $\Phi_{[s_1,s_2]_\d}=[[i_{s_1},\d],i_{s_2}(\s)]$. As $s \mapsto \Phi_s$  is injective, $[s_1,s_2]_\d$ is well defined. Moreover, using again the injectivity of $\Phi$  the Leibniz property for $[.,.]_\d$  is obtained from the following results:
$
\begin{array}{rl}
  \Phi_{[s_1,fs_2]_\d}(\s) & =\d(\s(s_1))(fs_2)-\d(\s(fs_2))(s_1)-\d\s(s_1,fs_2) \\
   & =f(\d(\s(s_1))(s_2)-\d(\s(s_2))(s_1)-\d\s(s_1,s_2))+\s(s_2)\d f(s_1) \\
   & =\Phi_{f[s_1,s_2]_\d}(\s)+\Phi_{\r(s_1)(f)}(\s) 
\end{array}
$

Now, if $\d^2=0$ by "formal argument" as in finite dimension, used for instance in \cite{Marl}, we can prove that $[.,.]_\d$  satisfies the Jacobi identity. From (\ref{rhodif}) we obtain that $\d=d_\r$ on $\cal F$. From Proposition \ref{crochetderiv}, we obtain that $d_\r=\d$ on $\G(\t_*)=\L^1\G^*(\t)$.
\end{proof}

\begin{rem}\label{pbderiv}${}$
\begin{enumerate}
\item Note  that, in general, if $M$ is not regular, a derivation  of the module of smooth functions $\cal F$ on a Banach manifold $M$ is not localizable (see for instance \cite{KrMi} section 35.1) . However, to our known there exists  no example of such a derivation which is not localizable. So, in Theorem \ref{algebroidetdiff}, if $M$ is not regular, we must impose that the  A-exterior differential on  $\L\G^*(\t)$ is  localizable.
\item In  Theorem \ref{algebroidetdiff}, we cannot assert that $\d$ and $d_\r$ coincides on $\L^k\G^*(\t)$ for $k\geq 2$.
Indeed, recall that  the Banach space $\L^k\E^*$ is generated by exterior products of $1$ forms  $\xi_{i_1}\wedge\cdots\wedge\xi_{i_k}$. However, even if $\E^*$ has a Schauder basis $\{\epsilon_\a\}_{\a\in \N}$ the family $$\{\epsilon_{\a_1}\wedge\cdots\wedge \epsilon_{\a_k},\;\; \a_1<\cdots<\a_k\}$$
{\bf could not be a Schauder basis} of $\L^k\E^*$
  for $k\geq 2$ (this is an unsolved  problem see \cite{Ram}).\\
\noindent  So,  for infinite dimensional  Banach spaces, locally, the module of local $k$-forms  $\L^k\G^*(\t)$ is not finitely generated but,  moreover  we have no "good topology" on $\L^k\G^*(\t)$ such that each $k$ form $\xi$ {\bf can not  be locally written as}
$$\eta(x)=\dis\sum_{i_1<\cdots< i_k}\eta_{i_1\cdots i_k}(x) \xi_{i_1}\wedge\cdots\wedge \xi_{i_k}$$ for some appropriate finite sequences of smooth functions $\xi_{i_1},\cdots \xi_{i_k}$
\noindent (for topologies on modules of sections see \cite{KrMi} and \cite{Lla})

 As a consequence,  any two A-derivations  which coincide on ${\cal F}(U)=\L^0\G^*(\t_U)$ and  on \\$\G^*(\t_U)=\L^1\G^*(\t_U)$ {\bf can be different} on $\L^k\G^*(\t_U)$, for $k\geq 2$. Once more, unfortunately, we have no example of such a situation.
 \end{enumerate}
\end{rem}

\subsection{Set of AL  structures on an anchored Banach bundle}
All the essential previous  results will be summarized  in the next theorem.

Let  $(E,\t,M,\r)$ be an anchored Banach bundle. We denote by  ${\cal ALB}(\t,\r)$ the set of  (localizable) AL-brackets on $(E,\t,M,\r)$ with fixed anchor morphism $\r$.


We have seen that to any sub AP-morphism on the dual bundle $E^*$ is associated an anchor morphism $\r_P:E\ap TM$ characterized  by
$$\textrm{ to } s \in \G(\t) \textrm{ one associates the derivation } f \mapsto \{\Phi_s,f\circ\t_*\} \textrm{ on } {\cal F}.$$
We denote by ${\cal AP}(\t,\r)$ the set of sub AP- Poisson morphisms on the dual bundle $E^*$ such that the associated anchor morphism is $\r$.

\begin{theo}${}$\\
Let $(E,\t,M,\r)$ be an anchored  Banach bundle. The set   ${\cal ALB}(\t,\r)$ has a natural structure of affine space in the following sense:
 
 given any $[.,.]_E\in {\cal ALB}(\t,\r)$ then we have:
 ${\cal ALB}(\t,\r)=\{ [.,.]_E+D,\;\; D\in \L^2\G(\t)\}.$\\
   
  \noindent  There exists a bijection  from   ${\cal ALB}(\t,\r)$ to ${\cal AP}(\t,\r)$ defined in the following way:\\
   at any  AL-bracket $[.,.]_E$ one associates the AL-bracket  $\{.,.\}_E$ on ${\cal F}^\f(E^*)$ characterized by $$\Phi_{[s_1,s_2]_E}=\{\Phi_{s_1},\Phi_{s_2}\}_E$$
  \end{theo}

  \noindent\begin{proof}\so{\it Proof}\\

  The only thing to prove is the structure of ${\cal ALB}(\t,\r)$. Consider two AL-brackets $[.,.]_i$, $i=1,2$, on $(E,\t,M,\r)$. It easy to see that  $D=[.,.]_1-[.,.]_2$  is an element of $\L^2\G(\t)$. The others properties come from
  Theorem \ref{linearAP}.\\
  \end{proof}

 \section{ Mechanical systems on an almost Lie algebroid}\label{mecha}
 \subsection{Hamiltonian system and Hamilton-Jacobi equation}

 Let $(E,\t,M,\r,[.,.]_\r)$ be an AL-algebroid. Denote by $P:T^\f E^*\ap TE^*$ the sub AP-morphism associated to this AL-strucuture and $\{.,.\}_P$ the associated  AP-bracket on ${\cal F}^\f(E^*)$ (see subsection \ref{APmor}). Any function $h\in {\cal F}^\f(E^*)$ is called a Hamiltonian function and the triple $(E, \{.,.\}_P,h)$  is called a {\bf Hamiltonian system}.  As we have already seen, to $h$ is associated a vector field $grad^P(h)=-P(dh)$ called the sub A-hamiltonian   gradient of $h$. As in this section, the AL-algebroid is fixed, $grad^P(h)$ will be denoted by $\overrightarrow{h}$. An integral curve of $\overrightarrow{h}$ is called a {\bf solution} of the Hamiltonian system  $(E, \{.,.\}_P,h)$\\

 \noindent\so{Local expressions} :\\

 In local trivialization (subsection \ref{loctriv}), given a Hamiltonian $h$, in local coordinates $(x,\xi)$ in $E_U$, we denote by $D_1h$ and $D_2h$ the partial derivative according to the variable  $x$ and $\xi$ respectively. The A-hamiltonian gradient of $\overrightarrow{h}$ has components $\overrightarrow{h}_1$ and $\overrightarrow{h}_2$, on $\M$ and $\E^*$ respectively. So we have the following characterization (see (\ref{Ploc})):
 \begin{eqnarray}\label{loctrivham}
 \begin{array}{c}
 <df',P_{\o(x)}(dh)>= <df', \r(D_{2}h)>\hfill{}\\
 <u',P_{\o(x)}(dh)>=<\xi,C_x(D_2h,u')>-<D_1h,R_x(u')>
 \end{array}
\end{eqnarray}

So we get
 \begin{eqnarray}\label{locham}
 \begin{array}{c}
 \overrightarrow{h}_1=\r(D_{2}h)\hfill{}\\
<u',\overrightarrow{h}_2>=<\xi,C(D_2h,u')>-<D_1h,R(u')>
\end{array}
\end{eqnarray}

 When $\M$ and $\E$ have unconditional basis, we have:
 $$\overrightarrow{h}=\dis\sum_{i\in I,\a\in A}\r^i_\a[\frac{\partial h}{\partial \xi_\a}\frac{\partial}{\partial x^i}-\frac{\partial h}{\partial x^i} \frac{\partial}{\partial \xi_\a} ]-\dis\sum_{\a,\b,\g\in A}C_{\a\b}^\g\xi_\a\frac{\partial h}{\partial \xi_\b}\frac{\partial}{\partial \xi_\a}$$
 So in local coordinates,  the integral curves of   $\overrightarrow{h}$ satisfies the following differential equations:
 \begin{eqnarray}
\dot{x}^i=\dis\sum_{\a\in A}\r_\a^i\frac{\partial h}{\partial \xi_\a} \qquad \dot{\xi}_\a=-\dis\sum_{i\in I}\r_\a^i\frac{\partial h}{\partial x^ i}-\sum_{\b,\g\in A}C_{\a\b}^\g\xi_\g\frac{\partial h}{\partial \xi_\b}
\end{eqnarray}
\bigskip
\bigskip

As in finite dimension (see for instance  \cite{LMM}) we have the following result on Hamilton-Jacobi equation

\begin{theo}\label{HJ}${}$\\
Let  $(E, \{.,.\}_P,h)$ be a Hamiltonian system. Given any section  $\o\in\G^*(\t)$	  we  denote by $\overrightarrow{h}_\o$ the vector field on $M$ defined by
$$\overrightarrow{h}_\o(x)=T_{\o(x)}\t_*(\overrightarrow{h}(\o(x))$$

Assume that $d_\r\o=0$; then the following properties are equivalent:
\begin{enumerate}
\item[(i)]  If $ c:I \ap M$ is an integral curve of the vector field $\overrightarrow{h}_\o$, then $\o\circ c:I\ap E^*$ is a solution of the Hamiltonian system  $(E, \{.,.\}_P,h)$.
\item[(ii)]  $\o$ satisfies the Hamilton-Jacobi equation
 i.e. $d_\r(h\circ \o)=0$
\end{enumerate}
 \end{theo}

 \noindent\begin{proof}\so{\it Proof of Theorem \ref{HJ}}  (adaption in our context of the proof of Theorem 4.1 of \cite{LMM})${}$\\

 We will  use the local trivializations (subsection \ref{loctriv}). So $\o(z)=(z,\xi(z))$ where $\xi$ is a smooth map from $U$ to $\E^*$. In this context,  ${\cal L}_\o(x)=T_x\o(\r(E_x))\subset T_{\o(x)}E^*\equiv \M\times \E^*$ is  the vector space $\{(\r(x,v),   D\xi(\r(x,v)), v\in \E\}\subset \M\times \E^*$. Denote by $[{\cal L}_\o(x)]^0\subset T^\f_{\o(x)} E^*\equiv \M^*\times \E$ the vector space
$$\{\eta=d(\Phi_s+f\circ\t_*(\o))\textrm{  such that } <\eta, \alpha>=0 \textrm { for all } \alpha \in {\cal L}_\o(x)\}.$$
In  our trivializations we have:
$$ [{\cal L}_\o(x)]^0=\{(df,v')\textrm{ such that } <df,\r(x,v)>+<D\xi(\r(x,v)),v'>=0 \textrm{ for all } v\in \E\}$$
\bigskip
\begin{lem}:\label{cal L}
\begin{enumerate}
\item $P([{\cal L}_\o(x)]^0)={\cal L}_\o(x)$ if and only if $d_\r \o=0$
\item if $d_\r \o=0$ then $\ker P_\o(x)\subset [{\cal L}_\o(x)]^0$
\end{enumerate}
\end{lem}

\noindent\begin{proof}\so{\it Proof of part 1}\\
In the previous trivializations we have
\begin{eqnarray}\label{calLoc}
<(df',v'), (\r(x,v),   D\xi(\r(x,v))>=<df',\r(x,v)>+< D\xi(\r(x,v)),v'>
\end{eqnarray}

So for any $(df,v)\in  [{\cal L}_\o(x)]^0$ according to (\ref{Ploc}) and (\ref{calLoc}), we have the following  equality for any $(df',v')\in T_{\o(x)}^\f E^*$
\begin{eqnarray}\label{calLoc=}
<df',\r(x,v)>+< D\xi(\r(x,v)),v'>=<\xi,[s,s']_\r>(x)+ <df', \r(x,v)>- <df, \r(x,v')>
\end{eqnarray}
 if and only if we have  $d_\r\o(s,s')(x)=0$, using the expression (\ref{ds}), for any $s'(z)=(z,v'(z))$.\\

The proof will be completed if we have (\ref{calLoc=}), for any given  local section  $ s:s(z)=(z,u(z))$, there exists a function $f_s$ such that  $(df_s,v)(x)$ belongs to $[{\cal L}_\o(x)]^0$.\\

Indeed, fix such a section $s$. We define  $\tilde{f}_s:\r(E_x)\subset T_xM\equiv \M\ap \R$  by
$$\tilde{f}_s(\r(x,u))=- <D\xi(\r(x,u)),v(x)>$$
for any $u\in E_x$. So $\tilde{f}_s$ is a linear form on $\r(E_x)$. From Hahn-Banach theorem, there exists on $\E$ a continuous linear form $f_s$ such that $f_s=\tilde{f}_s$ on $\r(E_x)$. So, it follows that we have

$<df_s,\r(x,u)>=- <D\xi(\r(x,u)),v(x)>$ for all $u\in E_x$ i.e. $(df_s,v(x))$ belongs to $[{\cal L}_\o(x)]^0$.\\

\noindent\so{\it Proof of part 2} : \\Consider $(df,v)\in \ker P(\o(x)) \subset T_{\o(x)}^\f E^*\equiv \M^*\times \E$ . So, for any $(df',v')\in T_{\o(x)}^\f E^*\equiv \M^*\times \E$, from (\ref{Ploc}) we have
\begin{eqnarray}\label{part2}
<\xi,[s,s']_\r>(x)+ <df', \r(x,v)>- <df, \r(x,v')>=0
\end{eqnarray}
Under the assumption $d_\r\o=0$ the relation  (\ref{part2}) is equivalent to:
\begin{eqnarray}\label{part3}
<D\xi(\r(x,v)),v'>-<D\xi(\r(x,v')),v>+ <df', \r(x,v)>- <df, \r(x,v')>=0
\end{eqnarray}
So for any given $(x,v')$ choose $f'$ such that $(df',v') $ belongs to $[{\cal L}_\o(x)]^0$. For this choice, we get

$<D\xi(\r(x,v'),v>+  <df, \r(x,v')>=0$.

\noindent It follows that $(df,v)$ belongs to $[{\cal L}_\o(x)]^0$\\
\end{proof}

\noindent {\it We come  back to the proof of Theorem \ref {HJ}}. The end of this proof follows exactly the same arguments as in Theorem 4.1 of \cite{LMM} with some adaptations. First property (i) is clearly equivalent to

(i') $T\o(\overrightarrow{h}_\o(x))=\overrightarrow{h}(\o(x))$

 We begin by the implication  (i')$\Rightarrow$ (ii).

 With the previous notations, $h$ is a function of the variables $(x,\xi)\in U\times \E^*$. Denote by $D_{2}h$  the partial derivative of $h$ according to variable $\xi$. As $dh$ belongs to $T^\f E^*$, the differential $D_2h$ gives rise to a section of $E_U$. The component of $\overrightarrow{h}_\o(x)=P_{\o(x})(dh)\in T_{\s(x)}^\f E^*\equiv \M\times \E^*$ on $\M$ is characterized by (see (\ref{locham})):

  $\overrightarrow{h}_\o(x)=\r(D_{2}h)(x)$ .

   \noindent From assumption (i'), we then have  $\overrightarrow{h}(\o(x))\in {\cal L}_\o(x)$.
From Lemma \ref{cal L} part 1, there exists $\eta\in [{\cal L}_\o(x)]^0$ such that
$$\overrightarrow{h}(\o(x))=P_{\o(x)}(\eta)$$
 So $(\eta-dh)(x)$ belongs to $\ker P_{\o(x)}$ and, using Lemma \ref{cal L} part 2  $(dh)(x)$ also belongs to  $[{\cal L}_\o(x)]^0$.
 But
 $$<d(h\circ\o),\r(s)>=<d h,T_x\o\circ\r(s)>$$
 for any $s\in E_x$.
 As $(dh)(x)$ belongs to $[{\cal L}_\o(x)]^0$ and $T_x\o\circ\r(s)$ belongs to ${\cal L}_\o(x)$ we get that $<d(h\circ\o),\r=0>$.\\

 (ii)$\Rightarrow$ (i').\\
 Under the assumption $d_\r\o=0$, if we have  $<d(h\circ\o),\r>=0$, as previously, we can show that $(dh)(x)$ belongs to $[{\cal L}_\o(x)]^0$, so $\overrightarrow{h}(\o(x))=P_{\o(x)}(dh)$, belongs to $P_{\o(x)}([{\cal L}_\o(x)]^0)$. But, using Lemma \ref{cal L} part 1, there exists $s\in E_x$ such that
 \begin{eqnarray}\label{hox}
\overrightarrow{h}(\o(x))=T_x\o(\r(s))
\end{eqnarray}
So we obtain
$$ \overrightarrow{h}_\o(x)=T_{\o(x)}\t_*(\overrightarrow{h}(\o(x))=T_x(\t_*\circ\o)(\r(s))=\r(s)$$
Using (\ref{hox}), we finally get
 $$T\o(\overrightarrow{h}_\o(x))=\overrightarrow{h}(\o(x))$$
  \end{proof}

\subsection{  Lagrangian  and Euler-Lagrange  equation on an AL-algebroid}\label{ALrie}

Given an anchored Banach bundle $(E,\t,M,\r)$, a {\bf semi spray} $S$ is a vector field on $E$ such that (see \cite{Ana}):

$p_E\circ S=Id_E$ where $p_E:TE\ap E$ is the canonical projection;

$T\t\circ S=\r$ where $T\t:TE\ap TM$ is the tangent map of $\t$.

On the other hand, a $C^k$-curve ($k\geq 1$) $c:[a,b]\ap E$ is called admissible, if we have $T\r( \dot{c}(t)=\r(c(t))$ for all $t\in [a,b]$. According to \cite{Ana}, we have:

\begin{prop}\label{adS}\cite{Ana}${}$\\
A vector field $S$ on $E$ is a semi spray if and only if all integral curves of $S$ are admissible curves.
\end{prop}

\noindent Among the class of semi sprays  the subclass of sprays  takes an important place for applications:
if we denote by $h_\l :E \ap E$, the homothety of factor $\l >0$ ($h_\l(u)=\l u$ for any $u\in E_x$ and any $x\in M$)
a  semispray $S$  is a spray if we have	$S\circ h_\l =\l Th_\l\circ S.$\\

\noindent\so{\it Local expressions}:\\

In local trivializations (subsection \ref{loctriv}), a semispray can be written in the  following way (see \cite{Ana}):
\begin{eqnarray}\label{locsemisp}
S(x,u)=(x,u, R_x(u),-2G(x,u))
\end{eqnarray}

When $\M$ and $\E$ have unconditional basis, we have:
\begin{eqnarray}\label{locsemispbase}
S=\dis\sum_{\a\in A}[(\sum_{i\in I}\r^\a_i\dis\frac{\partial}{\partial x^i})- 2G^\a\dis\frac{\partial}{\partial u^\a}]
\end{eqnarray}

\bigskip
\bigskip

A {\bf Lagrangian} on a Banach bundle $(E,\t,M)$ is a smooth map $L:E\ap \R$. We say that $L$ is homogenous of degree $k$ if we have:
$$L\circ h_\l=\l^kL.$$
The following result is classical (see for instance  \cite{AbMa} section 3.5)

\begin{lem}\label{Lagrang}${}$\\
Let $L$ be a Lagrangian, we denote by $L_x$ the restriction of $L$ to the fiber $E_x$. Then the map $\L_L:(x,u)\ap(x, dL_x(u))$  is a bundle morphism from $E$ to $E^*$
\end{lem}

We will say that $L$ is  {\bf regular } (resp. {\bf  strong regular}) (resp. {\bf  hyperregular}) if $\L_L$ is an injective  morphism (resp. isomorphism) (resp. diffeomorphism). Notice that when $L$ is strong regular then $\L_L$ is a local diffeomorphism. So $L$ is hyperregular if and only if $L$ is strong regular and  the restriction of $\L_L$ to each fiber is injective.

Denote by $h_\l :E \ap E$, the homothety of factor $\l >0$ ($h_\l(u)=\l u$ for any $u\in E_x$ and any $x\in M$).
As in finite dimension, let $\Theta $ be the Liouville field on $E$ which is the vector field whose flow is the homothety $\{h_\l\}_{\l\in \R}$. In {\it  local trivializations}
we have $\Theta(x,u)=(x,u,0, u)$ and when $\M$ and $\E$ have unconditional basis, we have $\Theta=\dis\sum_{\a\in A} u^\a\dis\frac{\partial}{\partial u^\a}$.\\

We denote by $H_L$ the  {\bf Lagrangian energy} associated to $L$ i.e.
$$H_L=dL(\Theta)-L$$

Given a regular Lagrangian $L$ on an AL-algebroid $(E,\t,M,\r,[.,.]_\r)$,  $\L_L$ is a local diffeomorphism. So for any $(x,u)\in E$, there exists an open neighborhood $U\times V \subset E_U$ of $(x,u)$ such that $(\L_L)_{| U\times V}$ is a diffeomorphism.\\

  {\bf Now suppose that $\E$ is reflexive}. Consider  a regular Lagrangian $L$ on $E$ and $U\times V\subset E$ an open set on which $\L_L$ is a diffeomorphism.  

Let be the function $ h_L=H_L\circ {\L_L}^{-1}$ on $\L_L(U\times V)$. Then $dh_L(\s)$ belongs to $T^*_\s E^*$ on $\L_L(U\times V)$
On $U\times V$ we can define the vector field  $\overrightarrow{L}$ characterized by
$$(\L_L)_*(\overrightarrow{H_L})=\overrightarrow{h_L}$$
and which is called the {\bf  local Euler-Lagrange vector field} of $L$ on $U\times V$. In particular, when $L$ is hyperregular, $\overrightarrow{L}$ is globally defined and called {\bf Euler-Lagrange vector field} of $L$.\\

\begin{theo}${}$
Consider a regular Lagrangian $L$ on $E$.
\begin{enumerate}
\item Any  curve $c=(\g,\m):[a,b]\ap U\times V$ is an integral curve of the local Euler-Lagrange vector field of $L$ if and only if it is a solution of the {\bf Euler-Lagrange} equations
\begin{eqnarray}\label{EuLa}
\dot{x}=R(u)\;\;\;\;\; \dis\frac{d}{dt}(D_2L)= R^t(D_1L)-D_2L(C(\;,u))
\end{eqnarray}
where $R^t:U\ap L(\M,\E^*)$ is the  field $x\ap (R_x)^t$ of transpose of $R_x\in L(\M,\E)$ and where  $C(\;,u)$ denote, for a fixed $u$,  the field of linear maps $x\ap [v\ap C_x(v,u)]$ (recall that $x \mapsto C_x$ is a field of bilinear maps)
\item If $L$ is hyperegular, the  Euler-Lagrange vector field $\overrightarrow{L}$  is a semi-spray. Moreover, if $L$ is homogenous of degree $2$ then $\overrightarrow{L}$ is a spray.\\
\end{enumerate}
\end{theo}

\begin{rem}\label{EuLabase}${}$\\
If $\M$ and $\E$ have unconditional basis, then in the associated coordinate systems, the Euler-Equation can be written in the following way:
\begin{eqnarray}\label{EuLab}
\dot{x}^i=\dis\sum_{\a\in A}\r^i_\a u^\a \;\;\;\;\; \dis\frac{d}{dt}(\dis\frac{\partial L}{\partial u^\a})=\dis\sum_{i\in I}\r^i_\a\dis\frac{\partial L}{\partial x^i}-\sum_{\b,\g\in A}C_{\a\b}^\g u^\b\dis\frac{\partial L}{\partial u^\g}
\end{eqnarray}
for any $i\in I$.
When $A$ and $I$ are finite sets of indexes, (\ref{EuLab}) is  the classical  Euler-Lagrange equation on the Lie algebroid (see for instance \cite{GMM})
\end{rem}

\noindent\begin{proof}\so{\it Proof}${}$\\

We again adopt the notations in  local trivializations (subsection \ref{loctriv}) . So $L$ is a function of variable $(x,u)$ and $\L_L$ is the map $(x,u)\ap (x, D_2L(x,u))$. For simplicity, we will denote this  map by $\L$ and for a fixed point $(x,u)\in E$ we denote by $(x,\xi)=\L(x,u)$
The tangent map $T\L$  of $\L$ is:
\begin{eqnarray}\label{TvL}
\begin{pmatrix}
Id&0\cr
D_{12}L&D_{22}L\\
\end{pmatrix}
\end{eqnarray}
So $[(T\L)^*]^{-1}$ is
\begin{eqnarray}\label{TvL*}
\begin{pmatrix}
Id&-D_{21}L\circ (D_{22}L)^{-1}\cr
0 &(D_{22}L)^{-1}\\
\end{pmatrix}
\end{eqnarray}

\noindent On the other hand, we have $H_L(x,u)=D_2L(x,u)(u)-L(x,u)$ so wet get
\begin{eqnarray}\label{HL}
dH_L(x,u)=(D_{12}L(x,u)(u,\;)-D_1L(x,u)(\;),D_{22}L(x,u)(u,\;))
\end{eqnarray}
So as $dh_L=[(T\L)^*]^{-1}\circ dH_L$, from (\ref{TvL*}) and (\ref{HL}) we get
\begin{eqnarray}\label{dhL}
D_1h_L(x,\xi)=-D_1L(x,\xi) \textrm{ and } D_2h_L(x,\xi)=u
\end{eqnarray}
From (\ref{loctrivham}) the  AL Hamiltonian $\overrightarrow{h_L}=([\overrightarrow{h_L}]_1,[\overrightarrow{h_L}]_2$ is characterized by
\begin{eqnarray}\label{hamhL}
[\overrightarrow{h_L}]_1(x,\xi)=R_x(u), \textrm{ and } ,[\overrightarrow{h_L}]_2=- D_2L(x,\xi)\circ C_x(u,\;)  + D_1L(x,\xi)\circ \r_x
\end{eqnarray}

Now as $\overrightarrow{H_L}=T\L(\overrightarrow{h_L})$ from (\ref{TvL}) and (\ref{hamhL}) we get

$[\overrightarrow{h_L}]_1(x,u)=R_x(u)$

 $D_{22}L([\overrightarrow{h_L}]_2(x,u))=-D_{12}L(x,u)\circ R_x(u)+D_1L(x,u)\circ R_x-D_2L(x,u)\circ C_x(\;,u)$

We can easily see that  these last equations  are equivalent to the Euler- Lagrange equations

\end{proof}

\subsection{Riemannian AL-algebroid and mechanical system }\label{riemal}
  Let $(E,\t,M)$ be a Banach bundle. Denote by $S^2T^*M\ap M$  the Banach bundle of symmetric bilinear form on $TM$. Recall that a global section $g$ of this bundle is called a  {\bf riemannian metric} on $E$, if  for for any $x\in M$, the bilinear form $g_x$ on $T_xM$ is positive definite,   i.e.  $g_x(u,u)> 0$ for any  $u\not=0$

To any riemannian metric $g$ on $E$, is associated a bundle morphism $g^\f:E\ap E^*$ defined by $g^\f(X)(Y)=g(X,Y)$. Of course, $g^\f$ is always injective. We say that $g$ is a {\bf strong}  riemannian metric if $g^\f$ is surjective. Note  that, in these conditions,  the Banach space $\E$ is isomorphic to a Hilbert space, and so {\bf $\E$ must be reflexive}.

We will say that the AL-algebroid  $(E,\t,M,\r,[.,.]_\r)$ is a {\bf riemannian} AL-algebroid if there exists a  riemannian metric $g$ on $E$.
 In this situation, as in finite dimension, we can consider the  Lagrangian system map $L:E\ap \R$ of a {\bf mechanical system } on $E$ given given by
 \begin{eqnarray}\label{mecasyst}
 L(s)=\dis\frac{1}{2}g(s,s)-V(\t(s))
\end{eqnarray}
 where  $\dis\frac{1}{2}g(s,s)$ is the "kinetic energy" and   $V:M\ap \R$ the "potential energy" of the mechanical system.  The associated Lagrangian   energy  is then:
 $$H_L(x,u)=\dis\frac{1}{2}g(s,s)+V(\t(s)).$$
 The Legendre transformation $\L_L$  is  $g^\f$. Of course, the Lagrangian $L$ is hyperregular, so the Lagrangian field $\vec{L}$ is well defined. Moreover  if $V\equiv 0$ then $\vec{L}$ is a spray.\\

 \noindent\so{\it Local expressions}\\

 \noindent In local trivialization (cf subsection \ref{loctriv}), we also denote by $x\ap g_x$ the field  of symmetric bilinear maps on $\E$ associated to $g$ and $x\ap g_x^\f: \E\ap \E^*$ the field of  associated isomorphisms. With these notations, we have
 $$L(x,u)=\dis\frac{1}{2}{g}_x(u,u)-V(x).$$
 According to the local expression (\ref{locsemisp})
 we can  write
 $$\vec{L}(x,u)=(x,u,R(u),-2G(x,u))$$  where
 $G$ is characterized  by:
 \begin{eqnarray}\label{ecritG}
g_x^\f(G(x,u))=\dis\frac{1}{2}\Bigl[<R_x^t\circ D_1{g}^\f_x(u),u>-\frac{1}{2}<R_x^t\circ D_1{g}_x(u,u),\;> -<R_x^t\circ dV,\;>-<{g}^\f_x(u),C_x.(\;,u);>\Bigr]
\end{eqnarray}

The Euler-Lagrange equation is given by
\begin{eqnarray}\label{Euler-Lagrange}
\left\{\begin{array} [c]{l}
  \dot{x}=R(u)\\
  \dot{u}=\dis\frac{1}{2}({g}^\f_x)^{-1}\Bigl[<R_x^t\circ D_1{g}^\f_x(u),u>-\frac{1}{2}<R_x^t\circ D_1{g}_x(u,u),\;> -<R_x^t\circ dV,.>-<{g}^\f_x(u),C_x(\;,u)>\Bigr]
  \end{array}
\right. 
\end{eqnarray}
In particular, when $\M$ and $\E$ have unconditional basis , the bilinear map can be written as a matrix ${g}=(g_{\a\b})_{\a,\b\in A}$ and we have according to  (\ref{ecritG}) can be written
 \begin{eqnarray}\label{ELinbasis}
\dis\sum_{\b\in A}g_{\a\b}G^\b=\dis\frac{1}{2}\sum_{\b,\g\in A}\Bigl [\sum_{i\in I}\frac{\partial g_{\a\b}}{\partial x^i}\r^i_\g-\frac{1}{2}\sum_{i\in I}\frac{\partial g_{\b\g}}{\partial x^i}\r^i_\a-\sum_{\d\in A}C^\d_{\a\b}g_{\d\g}\Bigr] u^\b u^\g -\sum_{i\in I}\frac{\partial V}{\partial x^i}\r^i_\a
 \end{eqnarray}
and, as in finite dimension,  the Euler Lagrange equations have an analogue  expression  which is left to the reader.\\

 Consider  a riemannian L-algebroid  $(E,\t,M,\r,[.,.]_\r)$ and $g$ its riemannian metric. If \\$\t':F\ap M$ is a Banach subbundle of $\t:E\ap M$, we can defined the complemeted Banach subbundle $F^\perp\ap M$ whose fiber is $F^\perp_x$ is the orthogonal in $E_x$ of $F_x$ relatively to the metric $g$. Let be  $\Pi :E\ap F$ is the natural morphism projection, we can define an AL-bracket on the set of section of $F$ by:
 $$[s_1,s_2]'=\Pi[s_1,s_2]_\r$$
 (see Example \ref{ex1} n$^0\;1$). So, if $\r'=\r_{|F}$,  for the induced metric $g'$ on $F$ induced by $g$, $(F,\t',M,\r',[.,.]')$ is a riemannian AL-algebroid.

\noindent  In this context,  denote by

 $i_F:F\ap E$ the canonical inclusion;

 $i_F^*:E^*\ap F^*$ the dual projection;

 $\Pi^*: E^*\ap F^*$ the dual morphism of $\Pi$;

 $P: T^*E^*\ap TE^*$ the P-morphism on $E$ associated to $[\;,\;]_\r$

 $P':T^*F^*\ap TF$ the AP-morphism associated  to $[\;,\;]'$.\\

 The Lagrangian  $L(x,u)=\dis\frac{1}{2}{g}_x(u,u)-V(x)$ on $E$ induces a Lagrangian $L'=L\circ i_F$ on $F$ which is a {\bf constrained} Lagrangian on $E$. As in finite dimensional case (see \cite{Marr}), we associate a  mechanical system on the AL-algebroid $(F,\t',M,\r',[.,.]')$  called   {\bf constrained  mechanical system} on $E$  which is obtain from the  {\bf  unconstrainted system associated} to $L$ on $E$ with the following relations :

 the Legendre transformation $\L_{L'}:TF\ap T*F$ satisfies  $\L_{L'}=i_F^*\circ \L_L\circ i_F$;

  The hamiltonian $h_{L'}=H_{L'}\circ (\L_{L'})^{-1}$ on $F^*$  is also  given by  $h_{L'}=h_L\circ \Pi^*$;

  The Lagrangian  vector field $\vec{L'}$ associated to $L'$ is also $\vec{L'}= T\Pi\circ \vec{L}\circ i_F$.

 \section{Constrained mechanical system and  Hilbert snakes}\label{exemple}
 \subsection{ The  context of an Hilbert snake}\label{cont}
 We will now present the problem of the Hilbert snake and apply the previous results on Riemannian-AL agebroid. The reader can find a complete description of this situation in \cite{PeSa}.

 In finite dimensional,  a  snake (of length L) is a (continuous) piecewise $C^1$-curve $S : [0,L] \ap  \R^d$, arc-length parameterized  so that  the origin $S(0) = 0\in \R^d$. According to \cite {Ro}, "charming  a snake" is a control problem so that   its  "head"   $S(L)$ describes   a given $C^1$-curve $c : [0,1]\ap \R^d$ in "minimal way" . More precisely
we look for a one parameter  family $\{S_t\}_{t\in [0,1]}$ such that $S_t(L)=c(t)$ for all $t\in [0,1]$ and such  that the family $\{S_t\}$ has an minimal infinitesimal  kinematic energy. This problem  has the precise following formulation:

Each snake $S$ of length $L$ in $\R^d$ can be given by a  piecewise $C^0$-curve $u : [0,L] \ap  \mathbb{S}^{d-1}$ so that $S(t)=\int_0^tu(\t)d\t$. We look for a $1$-parameter  family $\{u_t\}_{t\in [0,1]}$  so that the associated  family $S_t$ of snakes satisfies  $S_t(L)=c(t)$ for all $t\in [0,1]$ so that the infinitesimal kinematic energy $\dis\frac{1}{2}\int_0^L||\frac{d}{dt}u_t(s)||ds$ is minimal.\\

A generalization of this problem in the context of an separable Hilbert space is developed in \cite{PeSa}.  More precisely, given a separable Hilbert space $\field{H}$ we consider the smooth hypersurface  $\field{S}^\infty$ of element of norm $1$. As previously, an Hilbert snake of length $L$ is a continuous piecewise
$C^1$-curve $S : [0,L] \ap  \field{H}$, arc-length parameterized  so that  $S(0) = 0$.  Each such snake is again given  by a  piecewise $C^0$-curve $u : [0,L] \ap  \field{S}^\infty$ so that $S(t)=\int_0^tu(\t)d\t$.  Given a fixed partition $\cal P$ of $[0,L]$, the set  ${\cal C}^L_{\cal P}$ of such curves will be called the configuration set and carries a natural structure of Banach manifold: when ${\cal P}=\{0,L\}$, the set ${\cal C}^L_{\cal P}$ is an hypersurface of the Banach  ${\cal C}([0,L], \field{H})$ of continuous map from $[0,L]$ to $\field{H}$ with the classical norm $||\;||_{\infty}$; for the general case ${\cal P}=\{0=s_0,\cdots s_N=L\}$  then ${\cal C}^L_{\cal P}$ is canonically homemophic to the  product  $[{\cal C}([0,L], \field{S}^\infty]]^N$ and so we put on ${\cal C}^L_{\cal P}$ the corresponding Banach structure product.
 Notice that, on each tangent space $T_u{\cal C}^L_{\cal P}$ we also have an $L^2$ product :
 \begin{eqnarray}\label{L2}
<v,w>_{L^2}=\dis\int_0^L <v(s), w(s)>ds
\end{eqnarray}
Where $<\;,\;>$ is the inner product in $\field{H}$. Of course for the associated norm $||\;||_{L^2}$ , the normed space $(T_u{\cal C}^L_{\cal P}, ||\;||_{L^2})$ is not complete.\\

To any "configuration" $u\in {\cal C}^L_{\cal P}$ is naturally associated the "end map"
${\cal E}(u)=\dis\int_0^Lu(s)ds$. This map is smooth and its kernel has a canonical complemented subspace which is the orthogonal of $\ker T_u{\cal E}$ in $T_u{\cal C}^L_{\cal P}$ according to the inner  product   (\ref{L2}). We then get a closed distribution $\cal D$ on ${\cal C}^L_{\cal P}$.
As in finite dimension,  for a one parameter  family $\{u_t\}_{t\in [0,1]}$   the associated  family $S_t$ of snakes satisfies  $S_t(L)=c(t)$ for all $t\in [0,1]$ so that the infinitesimal kinematic energy $\dis\frac{1}{2}\int_0^L||\frac{d}{dt}u_t(s)||ds$ is minimal,  if $c(t)$ has a "lift "  $\tilde{c}$ in ${\cal C}^L_{\cal P}$ which is tangent to $\cal D$, called an "horizontal lift". So  the problem for the head of the Hilbert snake  to join an initial state $x_0$  to a final state $x_0$ can be transformed in the following "accessibiliy problem" :

Given a initial (resp; final) configuration $u_0$ (resp. $u_1$) in ${\cal C}^L_{\cal P}$, so that ${\cal E}(u_i)=x_i$, $i=1,2$, find a piecewise $C^1$ horizontal curve $\g:[0,T]\ap {\cal C}^L_{\cal P}$ (i. e. $\g$ is tangent to $\cal D$) and  which joins $u_0$ to $u_1$. \\

Given any configuration $u\in {\cal C}^L_{\cal P}$ we look for the accessibility set ${\cal A}(u)$ of  all configurations $v\in {\cal C}^L_{\cal P}$ which can be joined from $u$ by piecewise $C^1$ horizontal curves. In the context of finite dimension, in   \cite{Ro},  using arguments about  the   action of the Mo\"{e}bus group on ${\cal C}^L_{\cal P}$, it can be shown that  ${\cal A}(u)$  is the maximal integral manifold of a finite dimensional distribution on  ${\cal C}^L_{\cal P}$. Unfortunately, in the context of Hilbert space, the same argument does not work. Moreover, as we are in the context of infinite dimension for $\field{S}^\infty$, we cannot hope to get a finite dimensional distribution whose maximal integral manifolds is ${\cal A}(u)$. However, we can construct a canonical distribution $\bar{\cal D}$ modelled on Hilbert space, which is integrable and so that the accessibility set ${\cal A}(u)$ is a dense subset  of  the maximal integral manifold through $u$ of $\bar{\cal D}$ (\cite{PeSa}Theorem 4.1).  Moreover this distribution is minimal in some natural sense. In fact, when $\field{H}$ is finite dimensional, $\bar{\cal D}$ is exactly the finite distribution in \cite{Ro} whose leaves are the accessibility sets.\\

\subsection{AL - algebroid structure}\label{ALst}
To  any Hilbert basis $\{e_i,\;i\in \N\}$, we can associate   a family of global vector fields $\{E_i,i\in \N\}$ on  ${\cal C}^L_{\cal P}$ which generates ${\cal D}$ and the anounced  distribution $\bar{\cal D}$ is the Hilbert distribution generated by $\{E_i, [E_j,E_l], \; i,j,l\in \N\}$. As $\cal E$ is not a submersion everywhere,  it follows that $\cal D$ is not a subbundle of $T{\cal C}^L_{\cal P}$.  \\
On one hand if  $\L=\{(i,j)\in \N^2, i<j\}$,  we can consider the Hilbert space $\field{G}=l^2(\N)\oplus l^2(\L)$
Now, given any Hilbert basis $\{e_i,i\in \N\}$ of $\field{H}$ we define an anchored bundle $( {\cal C}^L_{\cal P}\times \field{G},\r, {\cal C}^L_{\cal P})$ by
$$\r(\s,\xi)=\dis\sum_{i\in \N} \s_iE_i+\sum_{(j,l)\in \L}\xi_{jl}[E_j,E_l]$$
Of course $\r$ is well defined surjective and do not depend of the choice of the Hilbert basis.\\

On the other hand, the Lie bracket of vector fields of the family $\{E_i, [E_j,E_l], \; i,j,l\in \N\}$ satisfies the following relations: (Lemma 4.3 \cite{PeSa} )

$[E_i,E_j](u)=<e_j,u>E_i(u)-<e_i,u>E_j(u)$ for any $u\in{\cal C}^L_{\cal P}$ and any $i,j\in \N$;

$[E_i[E_j,E_k]]=\d_{ij}E_k-\d_{ik}Ej$  for any $i,j,k\in \N$

  $[[E_i,E_j],[E_k,E_l]]=\d_{il}[E_j,E_k]+\d_{jk}[E_i,E_l]-\d_{ik}[E_j,E_l]-\d_{jl}[E_i,E_k]$ for any $i,j,k,l\in \N$.\\

 So, on   $\field{G}$ we define a  Lie algebra structure in the following way:

let be $(\epsilon_i)_{i\in \N}$ (resp. $(\epsilon_{ij})_{(i,j)\in \L}$ the canonical basis of $(l^2(\N)$ (resp. $(l^2(\L))$;

on $\field{G}$ we define a  Lie algebra structure in the following way:

let be $(\epsilon_i)_{i\in \N}$ (resp. $(\epsilon_{ij})_{(i,j)\in \L}$ the canonical basis of $(l^2(\N)$ (resp. $(l^2(\L))$;

according to the previous relations, we then define:

$[\e_i,\e_j]=\o_{ij}$, for all $i,j\in \N$

$[\e_i,\o_{jk}]=\d_{ij}\e_k-\d_{ik}\e_j$, for all $i\in \N$ and $(j,k)\in \L$

 $[\o_{ij},\o_{kl}]=\d_{il}\o_{jk}+\d_{jk}\o_{il}-\d_{ik}\o_{jl}-\d_{jl}\o_{ik}$, for all $(i,j) (kl)\in \L$.\\

\noindent For any $\s=\dis\sum \s_i\a_i$, $\s'=\sum \s'_j\e_j$ in $l^2(\N)$ and  $\xi=\sum \xi_{ij}\o_{ij}$,  $\eta=\sum\eta_{kl}\o_{kl}$ in $l^2(\L)$, naturally we can define:

$[\s,\s']=\dis\sum_{i,j\in \N}\s_i\s'_j[\e_i,\e_j]$

$[\s,\eta]=\dis\sum_{i\in \N, (k,l)\in \L}\s_i\eta_{kl}[\e_i,\o_{kl}]$

 $[\xi,\eta]=\dis\sum_{(i,j)\in \L,(k,l)\in \L}\xi_{ij}\eta_{kl}[\o_{ij},\o_{kl}]$.\\

Coming  back to the anchored bundle  $( {\cal C}^L_{\cal P}\times \field{G},\r, {\cal C}^L_{\cal P})$, each section $\varphi$  of the trivial bundle ${\cal C}^L_{\cal P}\times\field{G}\ap {\cal C}^L_{\cal P}$ can be identified with a map $\varphi:{\cal C}^L_{\cal P}\ap \field{G}$. So, on the set $\G(\field{G})$ of section of this trivial bundle we can defined a Lie bracket by:
$$[\varphi,\varphi'](u)=[\varphi(u),\varphi'(u)]+d\varphi(\r(u,\varphi'(u))-d\varphi'(\r{u,\varphi(u)})$$
 It follows that  $( {\cal C}^L_{\cal P} \times \field{G} ,\r,{\cal C}^L_{\cal P},[\;,\;])$ is a {\bf  Banach  Lie algebroid structure}  on ${\cal C}^L_{\cal P}$\\

In $\field{G}$ let  be $\pi:\field{G}\ap l^2(\N)$ the canonical projection  whose kernel is $l^2(\L)$ and denote again by $\pi:{\cal C}^L_{\cal P}\times \field{G}\ap {\cal C}^L_{\cal P}\times   l^2(\N)$ the associated  projection bundle. Again  any section of the trivial bundle ${\cal C}^L_{\cal P}\times l^2(\N)\ap{\cal C}^L_{\cal P} $ can be identified with a map from ${\cal C}^L_{\cal P}$ to $l^2(\N)$. Of course the set $\G(l^1(\N))$  of such sections is contained in $\G(\field{G})$. So, as in Example \ref{ex1} 1, on $\G(l^2(\N))$, we can define an almost Banach  Lie bracket by:
$$[[\varphi,\varphi']](u)=\pi([\varphi,\varphi'](u)).$$
So, if we denote by  $\theta$ the restriction of $\r$ to $l^2(\N)\times {\cal C}^L_{\cal P}$
we get an {\bf almost  Banach Lie algebroid structure} $({\cal C}^L_{\cal P}\times l^2(\N),\theta, {\cal C}^L_{\cal P},[\;,\;])$ on ${\cal C}^L_{\cal P}$. \\

On the other hand, the  distribution $\bar{\cal D}$ is a weak Hilbert integrable distribution, on ${\cal C}^L_{\cal P}$, this means that, for any $u\in {\cal C}^L_{\cal P}$, there exists an Hilbert manifold $N$ and a smooth injective map $f:N\ap  {\cal C}^L_{\cal P}$ such that (see \cite{Pel}):\\
$u$ belongs to $f(N)$,\\
  $T_xf: T_xN\ap T_{f(x)}{\cal C}^L_{\cal P}$ is injective ,\\
  $T_xf(T_xN)=\bar{\cal D}_{f(x)}$ for any $x\in N$\\
  {\it We say that $N$ is an integral manifold through $u$  }\\

  Given such  integral manifold   which maximal (for the inclusion),  the pull back $f_*\{ {\cal C}^L_{\cal P}\times \field{G}\}$ and  $f_*\{ {\cal C}^L_{\cal P}\times l^2(\N)\}$   can be identified with  $N\times \field{G}$ and $N\times l^2(\N)$  respectively; Then, $\r$ (resp. $\theta$) induces an anchor $\r_N:N\times \field{G}\ap TN$ (resp.  $\theta_N:N\times l^2(\N)\ap TN$). The barcket $[[\;,\;]]$ and(resp. the almost bracket   $[[\;,\;]]$) induces a bracket (resp.  an almost bracket) again denoted $[[\;,\;]]$.So we get  also a Banach Lie  algebroid $(N\times \field{G}, \r_N,N,[[\;;\;]])$ and an  almost  Banach Lie algebroid structure $(N\times l^2(\N),\theta_N,N,[[\;,\;]])$ on $N$. \\
  Now recall the following result of  \cite{PeSa} :

  \begin{prop} \label{PropN}${}$\\
 Let be $N$ a maximal integral manifold of $\bar{\cal D}$ and fix some  $u\in N$. Then we have the following properties
  \begin{enumerate}
  \item The set $\S({\cal E})$ at which $\cal E:{\cal C}^L_{\cal P}\ap \field{H}$ is not a submersion is a weak manifold of ${\cal C}^L_{\cal P}$ which is diffeomorphic to the projective space $\field{P}^\infty$ of  $\field{H}$. Its complementary ${\cal R}({\cal E})$ is an open dense set of ${\cal C}^L_{\cal P}$. Moreover, $N$ is a maximal integral manifold of $\bar{\cal D}$.
  \item[(1)] Assume that  $u\in \S({\cal E})$ then $N=\S({\cal E})$. Let be ${\cal L}_v$  the $1$-codimensional Hilbert  subspace  $[\ker \theta_N]_v^\perp \subset \{v\}\times l^2(\N)$ for any $v\in N$. Then ${\cal L}=\dis\cup _{v\in N}{\cal L}_v$ is a $1$-codimensional Hilbert  subbundle of $N\times l^2(\N)$ and the restriction $\psi_N$ of $ \theta_N$ to ${\cal L}$ is an isomorphism onto $TN$ and we have ${\cal D}_{| N}=TN$.
   \item [(2)] Assume that $u\in {\cal R}({\cal E})$. Then $N$ is contained in ${\cal R}({\cal E})$.
   Let be  $\field{V}_u$ the Hilbert subspace of $\field{H}$ generated by the set
$$\{u(t)-u(0),\; t\in [0,L]\}$$
 and choose an Hilbert   basis $\{e'_a,\; a\in A\}$ (resp $e'_b,\; b\in B\}$) of $[\field{V}_u]^\perp$ (resp. $\field{V}_u$). If $\L_u$ is the set of pair $(i,j)\in \L$ such that
 that $i$ or $j$ do not belongs to $A$, then $N$ is an Hilbert manifold modeled on $l^2(\N)\oplus l^2(\L_u)$  and is contained in $ {\cal R}({\cal E})$ .\\
  Let be ${\cal L}_v$ the orthogonal of $\ker [\r]_v\subset   \{v\}\times \field{G}$. Then ${\cal L}=\dis\cup _{v\in N}{\cal L}_v$ is a Hilbert  subbundle of $N\times \field{G}^2$   which contains $N\times l^2(\N)$ and the restriction of $\psi_N$ of $\Psi_N$ to ${\cal L}$ is an isomorphism on $TN$ \\
  Moreover, ${\cal L}$ contains $N\times l^2(N)$ and the restriction  of $\theta_N$ to $N\times l^2(N)$ is an isomorphism on ${\cal D}_{| N}$
 \end{enumerate}
 \end{prop}
 So, on $N$, if we denote by $[\;,\;]$ the usual Lie bracket and $p_N:TN\ap N$ the tangent bundle, we have {\it the  canonical  L algebroid} $(TN,p_N,N,Id_N,[\;,\;])$. When $u\in {\cal R}({\cal E}$,  on $N$, the distribution $F={\cal D}_{|N}$ is an hilbert subbundle $p':F\ap N$ of $p_N:TN\ap N$.
 When $u\in\S({\cal E})$ we have ${\cal D}_{| N}=TN$ and so again we can consider  $p':F={\cal D}_{| N}\ap N$ as an Hilbert subbundle of $TN$
\subsection{Constrained mechanical system}

    Let be $x$ and $y$ two states of the head of  the Hilbert snake which can be joined by an optimal curve (in the   sense of subsection \ref{cont}) and consider the set $\O(x,y)$ the set of optimal curves $c$ which joins $x$ to $y$.  Classically if $c\in \O(x,y)$ is defined on $[0,T]$ its kinematic energy is
$$E(c)=\dis\frac{1}{2}\int_0^T||\dot{c}(t)||^2 dt$$
Let be $E(x,y)=\inf_{c\in \O(x,y)}E(c)$. Assume that there exists $\bar{c}\in \O(x,y)$ such that  $E(\bar{c})=E(x,y)$, then $\bar{c}$ will be called an  {\it optimal minimizing curve} which joins $x$ to $y$.  Of course, such a curve is also an optimal minimizing curve between any pair of its points.  So we can look for the existence of optimal minimizing curve which begin at a given  original state $x$ of the head of the considered snake.\\
 On one hand, each  $c\in \O(x,y)$  has an horizontal lift $\g$ in ${\cal C}^L_{\cal P}$, and  $\g$ lies in $N$ As  ${\cal L}\ap N$ is an Hilbert subbundle of  $N\times  \field{G}\ap N$, we have a natural  riemannian metric $g$ on ${\cal L}$. From Proposition \ref{PropN}, the isomorphism $\psi_N$ gives rise to a riemannian  metric - again denoted by $g$ - on $TN$ and induces a riemannian metric $g'$ on $F={\cal D}_{| N}$. Note that the  inner product  induces by $g'$ on each fiber ${\cal D}_u$, $u\in N$, is exactly the inner product   induced  by $\r_u:{\cal D}_u\ap T_u\field{H}\equiv \field{H}$ of the canonical one on $\field{H}$. So we can defined - as in subsection \ref{ALrie} - an AL bracket $[\;,\;,]'$ on $F$ and we get an {\bf  AL algebroid} $(F, p',i_F, N,[\;,\;]')$\\

Now,  coming back to our problem of optimal minimizing curve on $\field{H}$.  To each  $c\in \O(x,y)$  has an horizontal lift $\g$ in ${\cal C}^L_{\cal P}$, and  $\g$ lies in $N$ for such a lift  $\g$ be an  of some optimal curve  $ c:[0,T]\ap \field{H}$, we have then
$$E(c)=\dis\frac{1}{2}\int_0^T||\dot{\g}||^2 dt.$$
 Finally, It follows that if $\bar{c}$ is   an optimal minimizing curve which joins $x$ and $y$, then  the associated lift  $\bar{\g}$   in $N$  is an extremal of the Lagrangian $L':F\ap \R$:
$$L'(v,\s)=\dis\frac{1}{2}||\s||^2$$
 on the { AL algebroid} $(F, p',i_F, N,[\;,\;]')$.\\

 So we get a {\bf constrained mechanical system} on the natural riemannian algebroid $(TN,p,N,Id_N,[\;,\;])$.
 It follows that such $\bar{\g}$ is a {\bf solution of the Euler-Lagrange equation} of $L'$. \\

{\it We will now give the differential system  satisfied by of such extremals in local coordinates.}\\

 Fix some $u\in {\cal R}({\cal E})$. According to Proposition  \ref{PropN}  $N$ is modelled on $l^2(\N)\oplus l^2(\L_u)$. So we have a local coordinates $ (\s,\xi)=\{(\s_i)_{i\in \n},(\xi_{jl})_{(j,l)\in \L_u}\}$.
We denote by $\{\dis\frac{\partial}{\partial \s_i},\frac{\partial}{\partial \xi_{jl}}\}$ the local Hilbert basis  of $TN$ associated to this coordinates system. With these notations, we  have:

$$E_i=\dis\frac{\partial}{\partial \s_i}+\sum_{(i,l)\in \L_u}\s_l\frac{\partial}{\partial \xi_{il}}-\sum_{(l,i)\in \L_u}\s_l\frac{\partial}{\partial \xi_{li}}$$

So from (\ref{Euler-Lagrange}) and (\ref{ELinbasis}),  the components $(\s_{i},\xi_{jl})$ of an extremal is given by

 \begin{eqnarray}
\left\{\begin{array} [c]{l}
  \ddot{\s}_i=0, \; i\in \N\\
  \ddot{\xi}_{jl}=\dot{\s}_{j}\dot{\s}_l,\; (j,l)\in \L_u\\
    \end{array}
\right.  \nonumber
\end{eqnarray}

Now, for $u\in \S({\cal E})$,  if we consider a basis $(e_i)_{i\in \N}$ of $\field{H}$ such that $u=\pm e_1$, then $N$ is modelled on $e_1^\perp$ (see Proposition \ref{PropN}). So, in local coordinates $\{\s_i,\; i>1\}$ the component of an extremal satisfies
$$ \ddot{\s}_i=0,\; i>1$$

\end{document}